\documentclass[12pt,fleqn]{article}

\usepackage{amssymb}
\usepackage{subfigure, graphicx}

\usepackage{latexsym}
\usepackage{graphicx}
\usepackage{amssymb, amsmath}
\usepackage{bm}
\usepackage{comment}

\usepackage{color}

\def\G{\Gamma}
\def\calP{{\cal P}}
\def\calS{{\cal S}}
\def\scF{{\cal F}}
\def\R{\mathbb{R}}
\def\Z{\mathbb{Z}}

\def\veca{\bm{a}}
\def\vecb{\bm{b}}

\def\vecx{\bm{x}}
\def\vecX{\bm{X}}

\def\vecT{\bm{T}}
\def\vecN{\bm{N}}

\def\Image{\mathop\mathrm{Image}}

\def\<{\langle}
\def\>{\rangle}
\def\ratio{\mbox{\sf r}}
\def\eps{\varepsilon}

\def\ds{\hbox{d}s}
\def\dS{\hbox{d}S}
\def\du{\hbox{d}u}

\newtheorem{lemma}{Lemma}
\newtheorem{theorem}{Theorem}
\newtheorem{proposition}{Proposition}
\newtheorem{remark}{Remark}

\title{Computational and qualitative aspects of motion of plane curves with a curvature adjusted tangential velocity}

\author{Daniel \v Sev\v{c}ovi\v{c}$^1$ \and   Shigetoshi Yazaki$^2$
}

\begin{document}

\maketitle
\footnotetext[1]{Faculty of Mathematics, Physics and Informatics, 
Comenius University, 842 48 Bratislava, Slovak Republic. 
{\tt sevcovic@fmph.uniba.sk}}
\footnotetext[2]{Faculty of Engineering, University of Miyazaki, 
1-1 Gakuen Kibanadai Nishi, Miyazaki 889-2192, Japan. 
{\tt yazaki@cc.miyazaki-u.ac.jp}}

\begin{abstract}

In this paper we investigate a time dependent family of plane closed Jordan curves evolving in the normal direction with a velocity which is assumed to be a function of the curvature, tangential angle and position vector of a curve. We follow the direct approach and analyze the  system of governing PDEs for relevant geometric quantities. We focus on a class of the so-called curvature adjusted tangential velocities for computation of the curvature driven flow of plane closed curves. Such a curvature adjusted tangential velocity depends on the modulus of the curvature and its curve average. Using the theory of abstract parabolic equations we prove local existence, uniqueness and continuation of classical solutions to the system of governing equations. We furthermore analyze geometric flows for which normal velocity may depend on global curve quantities like the length, enclosed area or total elastic energy of a curve. We also propose a stable numerical approximation scheme based on the flowing finite volume method. Several computational examples of various nonlocal geometric flows are also presented in this paper. 
\end{abstract}



\pagestyle{myheadings}
\thispagestyle{plain}
\markboth{
D. \v SEV\v{C}OVI\v{C} and  S. YAZAKI}{On a motion of plane curves with a curvature adjusted tangential velocity}

\section{ Introduction}
\label{sec-1}

In this paper we investigate a time dependent family of plane closed Jordan curves $\Gamma^t, t\in[0,T),$ evolving in the direction of the inner normal with a speed $v$, which  is assumed to be a function of the curvature $k$, tangential angle $\nu$ 
and position vector $\vecx\in\Gamma^t$, 
\begin{equation}\label{geomrov}
v= \beta(\vecx, k, \nu) \,. 
\end{equation}
Recall that the evolving family of plane curves having the normal velocity speed of the form (\ref{geomrov}) can be often found in various  applied problems, e.g. dynamics of phase boundaries in thermomechanics, material science (motion of bipolar loops \cite{Paus2009,Paus2010}), image and movie processing in computer vision theory (see e.g.  \cite{SrikrishnanCRS2007}). For a comprehensive overview of industrial applications of the geometric equation having  the form of (\ref{geomrov}) we refer to a book by Sethian~\cite{Se2}.

We analyze a system of governing PDEs for geometric quantities and propose a numerical method for computing the mean curvature flow of plane closed curves with nontrivial tangential redistribution of points along evolving curves. We focus on a class of so-called curvature adjusted tangential velocities for which the tangential velocity depends on the function of curvature and its average along the curve. Following \cite{BKPSTY}, the tangential speed is constructed as a linear combination of the asymptotic uniform tangential redistribution developed in \cite{MikulaS2004a} and the tangential velocity extracted from the crystalline curvature flow equations by the second author in \cite{Yazaki20XXa}. Using the theory of abstract parabolic equations due to Angenent \cite{A2,A4} we show local existence and uniqueness of a classical smooth solution to the system of governing equations. We furthermore propose a numerical approximation scheme based on the flowing finite volume method. We present several computational examples. In contrast to the simplified numerical algorithm proposed in \cite{SY2008} (see also \cite{BKPSTY}), the curvature and tangent angle are not calculated from the position vector but the parabolic equations for the curvature and tangent angle are solved separately. Such an approach yields a stable and robust numerical approximation scheme derived by Mikula and the first author in \cite{MikulaS2001,MikulaS2004a,MikulaS2004b} for the case of the so-called (asymptotically) uniform tangential redistribution. 
Recently, a higher order discretization scheme involving asymptotically uniform tangential redistribution has been proposed and analyzed by Mikula and Bala\v{zo}vjech in \cite{BM2011}.

An embedded closed plane Jordan curve $\G$ can be parameterized by a smooth function $\vecx:\ \R/\Z\supset[0, 1]\to\R^2$ such that 
$\G=\Image(\vecx)=\{\vecx(u);\ u\in[0, 1]\}$ and 
$|\partial_u\vecx|>0$. 
We denote $\partial_\xi{\sf F}=\partial{\sf F}/\partial\xi$, and 
$|\veca|=\sqrt{\veca.\veca}$ where 
$\veca.\vecb$ is the Euclidean inner product between vectors $\veca$ and $\vecb$. 
The unit tangent vector is $\vecT=\partial_u\vecx/|\partial_u\vecx|=\partial_s\vecx$, 
where $s$ is the arc-length parameter $\ds=|\partial_u\vecx|\du$, and 
the unit inward normal vector is uniquely determined through the relation $\det(\vecT, \vecN)=1$. A signed curvature in the direction $\vecN$ is denoted by $k$. i.e. $k=\det(\partial_s\vecx, \partial_s^2\vecx)$. Let $\nu$ be the angle of $\vecT$, i.e., $\vecT=(\cos\nu, \sin\nu)$ and 
$\vecN=(-\sin\nu, \cos\nu)$. 
The problem of evolution of plane curves can be formulated as follows: 
given an initial curve $\G^0=\Image(\vecx^0)$, find a $t$-parameterized family of plane curves  $\{\G^t\}_{t\geq 0}$,  $\G^t=\{\vecx(u, t);\ u\in[0, 1]\}$ starting 
from $\G^0$ and evolving according to the normal velocity
$v=\beta(\vecx, k, \nu)$. 
We follow the so-called direct approach in which we describe evolution of plane curves by an evolution equation for the position vector $\vecx$:\ 
$\partial_t\vecx=\beta\vecN+\alpha\vecT$, $\vecx(., 0)=\vecx^0(.)$. 
Here $\alpha$ is the tangential component of the velocity vector. 
Note that $\alpha$ has no effect on the shape of evolving closed curves, 
and the shape is determined by the value of the normal velocity $\beta$ only. 
Hence the simplest setting $\alpha\equiv 0$ can be chosen. 
Dziuk \cite{Dziuk1994} studied a numerical scheme for $\beta=k$ in this case. 
For general $\beta$, however, such a choice of $\alpha$ may lead to numerical instabilities caused by undesirable concentration of grid points. In order to obtain stable numerical computation, 
several nontrivial choices of $\alpha$ have been proposed.  In \cite{Kimura1994} Kimura proposed a redistribution scheme in the case $\beta=k$  by choosing tangential velocity $\alpha$ such that $|\partial_u\vecx(u,t)|=L^t$, where $L^t=\int_{\G^t}\,\ds$ is the total length of $\G^t$. Hou, Lowengrub and Shelley \cite{HouLS1994} derived the tangentail velocity (\ref{eq:tangential_velocity}) with $\varphi\equiv 1$ and $\omega\equiv 0$ (see Section 2). 
They also mentioned  (\ref{eq:tangential_velocity}) for general $\varphi$ with $\omega\equiv 0$. But they did not make explicit comments on the importance of redistribution of grid points. Another choice of tangential velocity $\alpha=-\partial_u(|\partial_u\vecx|^{-1})$ for the case 
$\beta=k$ has been proposed by Deckelnick \cite{Deckelnick1997}. In \cite{MikulaS2001} Mikula and the first author derived (\ref{eq:tangential_velocity}) with $\varphi\equiv 1$ and $\omega\equiv 0$ 
in general frame work of the so-called intrinsic heat equation for $\beta=\beta(k, \nu)$. 
Moreover, in \cite{MikulaS2004a} and \cite{MikulaS2004b}, they proposed a method of asymptotically uniform redistribution for $\beta=\beta(\vecx, k, \nu)$.  Besides aforementioned uniform distribution methods, under the so-called crystalline curvature flow,  grid points are not uniformly distributed  but their redistribution takes into account variations in the curvature. The second author extracted the tangential velocity $\alpha=-\partial_s\beta/k$ which is implicitly built in the crystalline curvature flow equation \cite{Yazaki20XXa}.  The asymptotically uniform redistribution is quite effective and valid for a wide range of applications.  However, from a numerical point of view,  there is no reason to take uniform redistribution automatically.  Hence a new way of curvature adjusted tangential redistribution can be considered in order to take into account the shape of evolved curves and variations in the modulus of the curvature. 

In this paper we furthermore  analyze qualitative and quantitative properties of nonlocal geometric flows in which the normal velocity $v=\beta$ is given by  
\begin{equation}\label{geomrovnonloc}
\beta(\vecx, k, \nu) \equiv \tilde\beta(\vecx, k, \nu) + \scF_\Gamma,
\end{equation}
where $\tilde\beta(\vecx, k, \nu)$ is a local part of the normal velocity of an evolving curve $\Gamma$ at a point $\vecx$ locally depending  on the curvature $k$ and tangent angle $\nu$ at position $\vecx$ and $\scF_\Gamma\in\R$ is a nonlocal part of the normal velocity at $\vecx$ depending on the entire shape of the curve $\Gamma$. Typically $\scF_\Gamma$ is a function depending on the total length $L$, the enclosed area $A$ and the total elastic energy ${\mathcal E}=\int_\Gamma k^2 \ds$ of a curve $\Gamma$, i.e. $\scF_\Gamma= \scF(L,A, {\mathcal E})$. 

The paper is organized as follows. In section~\ref{sec:GE} we investigate the system of governing partial differential equations for relevant geometric quantities. We derive a tangential velocity taking into account variations in the modulus of the curvature along evolving curves. In section~\ref{sec:local} we prove, locally in time, existence and uniqueness of classical smooth solutions.  In section~\ref{sec:affine} we discuss the relationship between crystalline curvature adjusted tangential velocity and affine scale space evolution of planar curves studied by Angenent, Sapiro and Tannenbaum in \cite{AST1998}. We show that the tangential velocity implicitly incorporated in the affine intrinsic heat equation is the same as the one present in the crystalline curvature motion. Section~\ref{sec:nonlocal} we discuss several interesting examples of geometric flows in which the normal velocity depends on the total length and the enclosed area. A special attention is put on analysis of the gradient flow for the isoperimetric ratio. A numerical approximation scheme based on the flowing finite volume method is presented in section~\ref{sec:scheme} together with computational examples (section~\ref{sec:results}). 

\section{System of governing equations and curvature adjusted tangential redistribution}
\label{sec:GE}

In what follows, we shall assume $\beta^\prime_k (\vecx, k, \nu)>0$. Then we can express the normal velocity $\beta$ in the form: 
\[
\beta=w(\vecx, k, \nu)k+F(\vecx, \nu),
\] 
where $w(\vecx, k, \nu)>0$ and $F(\vecx, \nu) = \beta(\vecx, 0, \nu)$.  By  $\beta^\prime_k$, $\beta^\prime_\nu$ and $\nabla_{\vecx}\beta$ we denote partial derivatives of $\beta$ with respect to $k$ and $\nu$, and,  $\nabla_{\vecx}\beta=(\partial_{x_1}\beta, \partial_{x_2}\beta)^T$ for $\vecx=(x_1, x_2)^T$. 

According to \cite{MikulaS2004a} (see also \cite{MikulaS2001, MikulaS2004b}), one can derive a system of PDEs governing evolution of plane curves satisfying  $\partial_t\vecx=\beta\vecN+\alpha\vecT$. It is easy to check the following facts. Indeed, it follows from the transformation $\partial_u{\sf F}=g\partial_s{\sf F}$ ($g=|\partial_u\vecx|$ is the so-called local length), Frenet's formulae $\partial_s\vecT=k\vecN$, $\partial_s\vecN=-k\vecT$, commutation relation $\partial_t \partial_s - \partial_s \partial_t  = (\beta k -\partial_s \alpha)\partial_s$, expression for the curvature $k=\det(\partial_s\vecx, \partial_s^2\vecx)$, tangent angle $(\cos\nu, \sin\nu)=\partial_s\vecx$ that 
the curvature $k$, tangent angle $\nu$, local length $g$ and position vector $\vecx$ satisfy the following system of PDEs:
\begin{eqnarray}
&& \partial_t k=\partial_s^2\beta+\alpha\partial_sk+k^2\beta, 
\label{eq:equation-k} \\
&& \partial_t\nu=\beta_k'\partial_s^2\nu+(\alpha+\beta_\nu')\partial_s\nu+\nabla_{\vecx}\beta.\vecT,
\label{eq:equation-nu} \\
&& \partial_tg=\left(-k\beta+\partial_s\alpha\right)g,
\label{eq:equation-g} \\
&& \partial_t\vecx=w\partial_s^2\vecx+\alpha\partial_s\vecx+F\vecN,
\label{eq:equation-x}
\end{eqnarray}
for $u\in[0, 1]$ and $t>0$. A solution $(k,\nu,g,\vecx)$ to (\ref{eq:equation-k})--(\ref{eq:equation-x}) is subject to the initial condition $k(.,0)=k_0(.), \nu(.,0)=\nu_0(.), g(.,0)=g_0(.), \vecx(.,0)=\vecx_0(.)$ corresponding to the initial curve $\Gamma^0=\Image(\vecx_0)$ and periodic boundary conditions for $g, k, \vecx$ for $u\in[0, 1]\subset\R/\Z$ and $\nu(1, t)=\nu(0, t)+2\pi$, i.e. $\nu(0,t)\equiv\nu(1,t)$ mod $2\pi$. 

In what follows, we shall specify the tangential velocity function $\alpha$ by taking into account the curvature variation with respect to the averaged function of the curvature. First,  for the total length $L^t=\int_{\Gamma^t} \,\ds = \int_0^1 g(u,t) \du$ of a curve $\Gamma^t$ evolving in the normal direction with the speed $\beta$ we have the following identity:
\begin{equation}
\frac{d}{dt} L^t + \int_{\Gamma^t} k\beta \,\ds = 0
\label{eq:length}
\end{equation}
(see e.g. \cite{MikulaS2001}). Recall that the so-called relative local length $\ratio(u, t)=g(u, t)/L^t$  plays an important role in controlling redistribution of grid points (cf. \cite{MikulaS2001}). Let $T>0$ be the maximal time of existence of a solution. If $\ratio(u, t)\to 1$ uniformly w.r. to  $u$ as $t\to T$,  then redistribution of grid points becomes asymptotically uniform. In general, the uniform redistribution can stabilize numerical computation,  and it has no effect on the shape of evolved closed plane curves. Neither uniform nor asymptotically uniform redistribution take into account variations in the curvature along evolved curves. Therefore, a natural question arises:  how to  redistribute grid points densely (sparsely) on sub-arcs where the modulus of the curvature is large (small) with respect to its averaged value. To answer this question, we shall define the so-called $\varphi$-adjusted  relative local length quantity
\begin{equation}
\ratio_{\varphi}(u, t)=\ratio(u, t) \frac{\varphi(k(u, t))}{\<\varphi(k(., t))\> }
\equiv  \frac{g(u, t)}{L^t} \frac{\varphi(k(u, t))}{\<\varphi(k(., t))\>}, 
\label{ratio-phi}
\end{equation}
where $\<{\sf F}(., t)\>= (1/L^t) \int_{\G^t}{\sf F}(s, t)\,\ds$ is the arc-length average of a quantity $\sf F$. Throughout the paper we will assume the function $\varphi(k)$ is smooth and strictly positive, i.e. $\varphi(k)>0$. As for the computational experiments discussed in section~\ref{sec:results} we shall consider $\varphi(k)=1-\eps+\eps\sqrt{1-\eps+\eps k^2}$ for $\eps\in[0, 1)$. 
Notice that $\varphi(k)\to 1$ if $\eps\to 0^+$, $\varphi(k)\to |k|$ if $\eps\to 1^-$. Other choices of $\varphi$ are also possible. For example $\varphi(k) = \sqrt{\varepsilon^2 + |k|^{2m}}, m>0$. In \cite{MikulaS2001, MikulaS2004a, MikulaS2004b} the constant function $\varphi(k)\equiv 1$ has been utilized in order to construct (asymptotically) uniform tangential redistribution. On the other hand, the function $\varphi(k)=|k|$ is implicitly built-in the crystalline curvature flow as it was pointed out in \cite{Yazaki20XXa}. The desired $\varphi$-adjusted asymptotic uniform  redistribution can be obtained by choosing $\alpha$ in such a way that $\theta(u, t):=\ln\ratio_{\varphi}(u, t)\to 0$ as $t\to T$. To this end, let us suppose that $\theta(u,t)$ satisfies the equation:
\begin{equation}
\theta(u, t)=\ln((e^{\theta(u, 0)}-1)e^{-\int_0^t\omega(\tau)\,d\tau}+1).
\label{eq:theta}
\end{equation}
According to \cite{MikulaS2004a, MikulaS2004b} the shape function $\omega$ can be constructed in the form: $\omega =\kappa_1+\kappa_2\<k\beta\>,$ where $\kappa_1, \kappa_2\geq 0$ are given constants. With regard to (\ref{eq:length}) we have $\omega =\kappa_1-\kappa_2\partial_t\ln L^t$. In the case when the time of existence $T<\infty$ is finite and evolved curves shrinks to a point, i.e. $L^t\to 0$ as $t\to T$,  we can take $\kappa_1=0$ and $\kappa_2>0$. On the other hand, if $T=\infty$  we can take $\kappa_1>0$ and $\kappa_2=0$. In both cases we have $\int_0^T \omega(t) dt = +\infty$ and the convergence $\lim_{t\to T}\theta(u, t)=0$ is guaranteed.  Setting $\kappa_1=\kappa_2=0$, we conclude $\theta(u,t)$ is constant w.r. to the time $t$ yielding thus uniform tangential redistribution (cf. \cite{MikulaS2004a}).

Taking into account governing equations (\ref{eq:equation-k})--(\ref{eq:equation-x}), definition of the $\varphi$-adjusted  relative local length (\ref{ratio-phi}) we obtain, after straightforward calculations, that $\theta(u, t):=\ln\ratio_{\varphi}(u, t)$ satisfies equation (\ref{eq:theta}) iff 
the tangential velocity $\alpha$ satisfies: 
\begin{equation}\label{eq:tangential_velocity}
\frac{\partial_s(\varphi\alpha)}{\varphi}
=\frac{f}{\varphi}-\frac{\<f\>}{\<\varphi\>}
+\omega(t)\left(\ratio_\varphi^{-1}-1\right), \quad 
f=\varphi k\beta-(\partial_s^2\beta+k^2\beta)\varphi^\prime_k, 
\end{equation}
where $\varphi=\varphi(k)$ and $\varphi^\prime_k=\partial_k\varphi(k)$. To construct a unique solution $\alpha$, we assume the renormalization condition $\<\varphi(k) \alpha\>=0$, i.e. $\int_\Gamma \varphi(k) \alpha \ds =0$.

\section{Local existence and uniqueness of classical solutions}\label{sec:local}

In this section, by following ideas adopted from \cite{MikulaS2004b} (see also \cite{MikulaS2001,MikulaS2004a}), we shall prove a local time existence, uniqueness and continuation of a classical smooth solution to the governing system of equations (\ref{eq:equation-k})--(\ref{eq:equation-x}). However, we have to rewrite the system (\ref{eq:equation-k})--(\ref{eq:equation-x}) into its equivalent form which does not explicitly contains the derivative $\partial_s \alpha$ of the tangential velocity appearing in equation (\ref{eq:equation-g}) for the local length $g$. Its presence in (\ref{eq:equation-g}) is a technical  obstacle in a direct application of the general result \cite[Theorem 5.1]{MikulaS2004b} on existence and uniqueness of solutions to (\ref{eq:equation-k})--(\ref{eq:equation-x}). To this end, we note that the auxiliary function $\theta=\theta(u,t)$ given by equation (\ref{eq:theta}) is a solution to the ODE:
\[
\frac{d\theta}{dt} = (1-e^{-\theta}) \omega(t),
\]
where $\omega = \kappa_1 +\kappa_2 \<k\beta\>$. Therefore the function $\ratio_{\varphi}= e^\theta$ satisfies the ODE: $\partial_t \ratio_{\varphi} = (\ratio_{\varphi} -1) \omega(t)$.  With this replacement we obtain the following system of PDEs:
\begin{eqnarray}
&& \partial_tk=\partial_s^2\beta+\alpha\partial_sk+k^2\beta, 
\label{eq:equation-k2} \\
&& \partial_t\nu=\beta_k'\partial_s^2\nu+(\alpha+\beta_\nu')\partial_s\nu+\nabla_{\vecx}\beta.\vecT,
\label{eq:equation-nu2}
\\
&& \partial_t \ratio_{\varphi} = (\ratio_{\varphi} -1) (\kappa_1+\kappa_2\<k\beta\>),
\label{eq:equation-g2} 
\\
&& \partial_t\vecx=w\partial_s^2\vecx+\alpha\partial_s\vecx+F\vecN,
\label{eq:equation-x2}
\end{eqnarray}
which is equivalent to the system (\ref{eq:equation-k})--(\ref{eq:equation-x}). Indeed, since
\begin{equation}
\ratio_{\varphi}=\ratio \frac{\varphi(k)}{\<\varphi(k)\>}
=\frac{g}{L} \frac{\varphi(k)}{\<\varphi(k)\>},
\label{g-ratio}
\end{equation}
we can express differentials $\partial_s$ and $\,\ds$ in terms of $\partial_u$ and $\du$ as follows:
\begin{equation}
\frac{\partial}{\partial s} 
= \frac{1}{g}\frac{\partial}{\partial u}
= \frac{\varphi(k)}{L \<\varphi(k)\>}\frac{1}{\ratio_{\varphi}}\frac{\partial}{\partial u}
\quad\hbox{and}\quad
\,\ds = \frac{L \<\varphi(k)\>}{\varphi(k)} \ratio_{\varphi} \du.
\label{partials-gr}
\end{equation}

A solution $\Phi=(k, \nu, \ratio_{\varphi} , \vecx)$ to the system of PDEs (\ref{eq:equation-k2})--(\ref{eq:equation-x2}) is subject to the initial condition $\Phi(.,0) = \Phi_0$ corresponding to the initial curve $\Gamma_0=\Image(\vecx_0)$. 

In what follows, we recall key ideas of the abstract theory of nonlinear analytic semigroups developed by Angenent \cite[Theorem 2.7]{A4}. By means of this theory we will be able to prove local existence and uniqueness of smooth solutions to the system of governing equations (\ref{eq:equation-k2})--(\ref{eq:equation-x2}) which can be rewritten as an abstract nonlinear equation of the form 
\begin{equation}
\partial_t \Phi = {\mathcal G}(\Phi),\quad \Phi(0)=\Phi_0\,,
\label{abstrakteq}
\end{equation}
where ${\mathcal G}(\Phi) = G(\Phi,\alpha(\Phi))$ and $G(\Phi,\alpha)$ is the right hand 
side of (\ref{eq:equation-k2})--(\ref{eq:equation-x2}) with $\alpha=\alpha(\Phi)$ being the unique solution to equation (\ref{eq:tangential_velocity}) satisfying the renormalization condition $\<\varphi(k) \alpha\>=0$. 

Now, suppose that the Fr\'echet derivative ${\mathcal G}^\prime(\bar\Phi)\in{\mathcal L}(E_1,E_0)$ belongs to the so-called maximal regularity class ${\mathcal M}_1(E_0,E_1)$ for any 
$\bar\Phi\in{\mathcal O}_1$ where ${\mathcal O}_1$ is an open neighborhood of the
initial condition $\Phi_0$ (cf. \cite{A4}). Here $E_0, E_1$ are Banach spaces, $E_1$ is densely embedded in $E_0$. Then, it follows from \cite[Theorem 2.7]{A4} that the abstract initial value problem (\ref{abstrakteq}) has a unique solution  $\Phi\in Y^{(T)} \equiv C([0,T], E_1) \cap  C^1([0,T], E_0)$ on some small enough time interval $[0,T]$. A maximal regularity class ${\mathcal M}_1(E_0,E_1) \subset {\mathcal L}(E_1,E_0)$ consists of those generators of analytic semigroups $A: D(A)=E_1 \subset E_0 \to E_0$ for which the linear equation $\partial_t \Phi = A\Phi + h(t)$,  $t\in(0,1]$, $\Phi(0)=\Phi_0$, has a unique solution $\Phi\in Y^{(1)}$, for any $h\in
C([0,1], E_0)$ and $\Phi_0\in E_1$. The pair $(E_0,E_1)$ is the so-called {\it maximal regularity pair}. It was also pointed out by Angenent in \cite{A4} (see also \cite{A2}), that the pair $(c^\varrho, c^{2+\varrho})$ of the so-called little H\"older spaces is the maximal regularity pair for the second order differential operator $A=-\partial^2_u$ subject to periodic boundary conditions at $u=0,1$. The ``little" H\"older space $c^{2 \mu+\varrho}=c^{2 \mu + \varrho} (S^1)$ with $0<\varrho<1$ and $\mu=0,\, 1/2,\, 1,$ is the closure of $C^\infty(S^1)$ in the topology of the H\"older space  $C^{2 \mu+\varrho}(S^1)$ (see \cite{A2}). It means $A\in{\mathcal M}_1(c^\varrho, c^{2+\varrho})$. 
In order to apply the abstract existence result \cite[Theorem 2.7]{A4} we define the following scale of Banach space
\[
E_\mu = c^{2 \mu + \varrho} \times c^{2 \mu + \varrho}_*  \times c^{1+\varrho}  \times (c^{2 \mu+\varrho})^2
\quad \hbox{for}\ \mu=0,\, 1/2,\, 1.
\]
By $c^{2 \mu + \varrho}_* (S^1)$ we denoted the Banach manifold 
$c^{2 \mu + \varrho}_* (S^1) = \{\nu:\R\to\R\,, \vecT =(\cos\nu,\sin\nu)
\in (c^{2 \mu + \varrho}(S^1))^2\}$. In other words, it is the space of all tangent angles $\nu\in c^{2 \mu + \varrho}$ satisfying the periodic boundary condition $\nu(0)\equiv \nu(1) \  \hbox{mod}\ 2\pi $.
 
Henceforth we shall assume $\varphi=\varphi(k)$ and $\beta=\beta(\vecx,k,\nu)$ are at least $C^3$ smooth functions such that $\varphi(k)>0$ and $\beta$ is a $2\pi$-periodic function in the $\nu$ variable. Furthermore, by ${\mathcal O}_{\frac{1}{2}}\subset E_{\frac{1}{2}}$ we shall denote a bounded open subset such that $\ratio_{\varphi}>0$ for any $(k,\nu,\ratio_{\varphi},\vecx)\in {\mathcal O}_{\frac{1}{2}}$. 

\begin{lemma} 
\label{lemma1}
Let $\alpha=\alpha(\Phi)$ be the tangential velocity function given as a unique solution to  (\ref{eq:tangential_velocity}) satisfying the renormalization condition $\<\varphi(k) \alpha\>=0$ where  $\Phi=(k, \nu, \ratio_{\varphi} ,x) \in {\mathcal O}_{\frac{1}{2}}$. Then 
$\alpha\in C^1(  {\mathcal O}_{\frac{1}{2}}, c^\varrho(S^1))$.

\end{lemma}

\noindent P r o o f: 
The term $f=\varphi k\beta-(\partial_s^2\beta+k^2\beta)\varphi^\prime_k$ appearing in the definition (\ref{eq:tangential_velocity}) of $\alpha$ can be decomposed as follows:
\begin{equation}
f=f_0 -\partial_s f_1,\quad \hbox{where} \ \ 
f_0= \varphi k\beta - k^2\beta \varphi^\prime_k + (\partial_s\beta)(\partial_s k) \varphi^{\prime\prime}_{kk}, \quad
f_1 = (\partial_s\beta) \varphi^{\prime}_{k}.
\label{deff0f1}
\end{equation}
Assuming $\varphi(k)>0$ and $\beta=\beta(\vecx,k,\nu)$ are $C^3$ smooth functions, we conclude $f_0=f_0(\Phi), f_1=f_1(\Phi) \in  C^1(  {\mathcal O}_{\frac{1}{2}}, c^\varrho(S^1))$. Since any averaged quantity (e.g. $\<\varphi(k)\>$, $\<\varphi k \beta\>$, $\<k \beta\>, \<f\>$) is constant in the arc-length variable $s$ it belongs to the class $C^1(  {\mathcal O}_{\frac{1}{2}}, c^\varrho(S^1))$. From (\ref{eq:tangential_velocity}) we have 
\begin{equation}
\partial_s(\varphi\alpha + f_1) 
= f_0 - \frac{\<f\>}{\<\varphi\>}\varphi 
+ \left(\ratio_\varphi^{-1}-1\right) \omega \varphi.
\label{alphaf0f1}
\end{equation}
Hence $\varphi\alpha + f_1 \in C^1(  {\mathcal O}_{\frac{1}{2}}, c^{1+\varrho}(S^1))$. Since 
$\varphi(k)>0$ and $f_1\in  C^1(  {\mathcal O}_{\frac{1}{2}}, c^\varrho(S^1))$ we conclude $\alpha\in C^1(  {\mathcal O}_{\frac{1}{2}}, c^\varrho (S^1))$, as claimed. 
\hfill$\diamondsuit$

\begin{remark}
\rm In \cite{MikulaS2001,MikulaS2004a} the authors proved that for the (asymptotically) uniform tangential velocity we have $\alpha\in C^1(  {\mathcal O}_{\frac{1}{2}}, c^{2+\varrho}(S^1))$, i.e. $\alpha(\Phi)\in c^{2+\varrho}(S^1)$ for $\varphi\in {\mathcal O}_{\frac{1}{2}}\subset E_{\frac{1}{2}}$. Indeed, in this case $\varphi$ is constant and therefore $f_1\equiv0, f_0\equiv k\beta \in C^1({\mathcal O}_{\frac{1}{2}}, c^{1+\varrho}(S^1))$. Such a higher smoothness was needed because the authors applied the abstract result \cite[Theorem 2.7]{A4} to the original system of equations (\ref{eq:equation-k})--(\ref{eq:equation-x}). Hence, in order to treat equation (\ref{eq:equation-g}) for  $g\in c^{1+\varrho}$ the tangential velocity $\alpha$ should belong to the class $c^{2+\varrho}(S^1)$.
\end{remark}

\medskip
Now we are in a position to state the main result of this section regarding existence, uniqueness and continuation of a smooth solution to the initial value problem (\ref{eq:equation-k2})--(\ref{eq:equation-x2}). 

\medskip
\begin{theorem}\label{localexistence}
Assume $\Phi_0=(k_0,\nu_0, \ratio_{\varphi0},\vecx_0)\in E_1$ where $k_0$ is the curvature, $\nu_0$ is the tangential vector, $\ratio_{\varphi0}>0$ is the $\varphi$-adjusted relative local length of an initial regular curve $\Gamma^0=\Image(\vecx_0)$. Assume $\varphi(k)>0$ and $\beta\equiv \tilde\beta(\vecx,k,\nu) + \scF_\Gamma$ where $\tilde\beta:\R^2\times\R\times\R\to\R$ and $\varphi:\R\to\R$ are $C^3$ smooth functions of their arguments such that $\tilde\beta$ is a $2\pi$-periodic function in the $\nu$ variable and $\min_{\Gamma_0}\tilde\beta^\prime_k(\vecx_0,k_0,\nu_0) >0$. The nonlocal part of the normal velocity $\scF_\Gamma$ is assumed to be a $C^1$ smooth function from a neighborhood ${\mathcal O}_{\frac{1}{2}}\subset E_{\frac{1}{2}}$ of $\Phi_0$ into $\R$, i.e. $\scF_\Gamma\in C^1({\mathcal O}_{\frac{1}{2}} ,\R )$.

Then there exists a unique  solution $\Phi=(k,\nu, \ratio_\varphi, \vecx) \in C([0,T], E_1) \cap  C^1([0,T], E_0)$ of the governing system of equations (\ref{eq:equation-k2})--(\ref{eq:equation-x2}) defined on some time interval $[0,T]\,,\ T>0$. 

\end{theorem}

\medskip
\noindent P r o o f: 
Since $\scF_\Gamma$ is a real valued nonlocal functional, $\scF_\Gamma: {\mathcal O}_{\frac{1}{2}} \to \R$,  we have $\partial_s\scF_\Gamma = 0$. Thus $\partial_s\beta =\partial_s\tilde\beta$. Similarly as in the proof of \cite[Theorem 5.1]{MikulaS2004a} we first rewrite equation (\ref{eq:equation-k2}) for the curvature in the form  
\begin{equation}
\partial_t k  =  \partial_s(\tilde\beta_k^\prime \partial_s k)
+ \partial_s(\tilde\beta_\nu^\prime k) + \partial_s( \nabla_{\vecx}\tilde\beta . \vecT)
+ \alpha \partial_s k + k^2\beta\,.
\label{curvatureeq}
\end{equation}
Here we have used the relations $\partial_s\nu = k$ and $\partial_s\tilde\beta = \tilde\beta_k^\prime \partial_s k +\tilde\beta_\nu^\prime k  +\nabla_{\vecx}\tilde\beta . \vecT$,  where $\vecT =(\cos\nu,\sin\nu), \vecN=(-\sin\nu,\cos\nu)$. Furthermore, using Frenet formulae $k\vecN =\partial_s \vecT =\partial_s^2 \vecx$ we can rewrite the position vector equation as follows: 
\[
\partial_t \vecx = \partial_s^2 \vecx + c \vecN +\alpha \vecT,
\] 
where $c(\vecx,k,\nu)\equiv \beta(\vecx,k,\nu) -  k$. This ``cheap trick" enables us to utilize the parabolic smoothing effect of the parabolic equation for $\vecx$ when we compared with an argument based only on the analysis of solutions to the ODE: $\partial_t \vecx = \beta \vecN +\alpha \vecT$ having not enough smoothness such as $\alpha(.) \in c^\varrho(S^1)$.

There exists an open neighborhood ${\mathcal O}_1\subset E_1$ of $\Phi_0$ such that 
$\Phi_0 \in {\mathcal O}_1$, and  $\ratio_\varphi>0$, $\tilde\beta^\prime_k(\vecx,k,\nu) >0$ for any $(k,\nu,\ratio_\varphi,\vecx)\in  {\mathcal O}_1$.  Then the mapping ${\mathcal G}$ is a $C^1$ smooth mapping from ${\mathcal O}_1 \subset E_1$  into $E_0$. Its linearization ${\mathcal G}^\prime$ at  $\bar\Phi =  (\bar k,\bar\nu, \bar\ratio_\varphi, \bar\vecx ) \in {\mathcal O}_1$ has the form:
\[
A={\mathcal G}^\prime(\bar\Phi) =A_1 +A_2, \quad \hbox{where}\ \  A_1= \bar D\, \partial^2_u, \quad A_2 = \bar B\, \partial_u + \bar C. 
\]
Here $\bar D$ is a diagonal matrix, $\bar D=\hbox{diag}(\bar D_{11},\bar D_{22}, \bar D_{33},\bar D_{44},\bar D_{55})$. The coefficients of $\bar D$ satisfy: $\bar D_{ii}\in c^{1+\varrho}(S^1)$,\ $\bar D_{11} = \bar D_{22} = \bar g^{-2}\tilde\beta^\prime_k(\bar \vecx,\bar k,\bar\nu)$, $\bar D_{33}=0$, $\bar D_{44} = \bar D_{55} = \bar g^{-2}\in c^{1+\varrho}(S^1)$, where $\bar g$ is given by (\ref{g-ratio}), i.e.  $\bar g = \bar\ratio_{\varphi} \bar L \<\varphi(\bar k)\>/\varphi(\bar k)$. The coefficients $\bar B_{ij}, \bar C_{ij}$ of $5\times5$ matrices $\bar B, \bar C$ corespond to the first and zero-th order derivative terms in the second order linear differential operator $A={\mathcal G}^\prime(\bar\Phi)$. 

Let ${\mathcal O}_{\frac{1}{2}}\subset E_{\frac{1}{2}}$ be an open subset in $E_{\frac{1}{2}}$ such that ${\mathcal O}_1= {\mathcal O}_{\frac{1}{2}}\cap E_1 \subset E_1$. With regard to Lemma~\ref{lemma1} we have $\alpha\in C^1({\mathcal O}_{\frac{1}{2}}, c^{\varrho} (S^1))$. Furthermore, we assumed $\scF_\Gamma\in C^1({\mathcal O}_{\frac{1}{2}} ,\R )$. As a consequence we obtain $\bar B_{ij}, \bar C_{ij}\in c^\varrho(S^1)$ for $i,j=1,\cdots,5$. 

For the second order differential linear operator $A_1$ we have $A_1\in {\mathcal L}(E_1, E_0)$.  Moreover, $A_1$, is a generator of an analytic semigroup on $E_0$ with the domain  $D(A_1)=E_1\subset E_0$. It belongs to the maximal regularity pair $(E_0,E_1)$, i.e $A_1\in {\mathcal M}_1(E_0,E_1)$ (see \cite{A4}). 

Since $A_2$ contains differentials of the first order only, it can be extended to a bounded operator from $E_{\frac{1}{2}} \supset E_1$ into $E_0$. Due to boundedness of $A_2\in {\mathcal L}(E_{\frac{1}{2}}, E_0)$ and taking into account the interpolation inequality between $c^\varrho$ and $c^{2+\varrho}$ spaces we conclude the following inequality:
\[
\Vert A_2 \Phi\Vert_{E_0}\le C_0 \Vert\Phi\Vert_{E_{\frac{1}{2}}}  
\le C_0 \Vert\Phi\Vert_{E_1}^{1/2}  \Vert\Phi\Vert_{E_0}^{1/2},
\] 
where $C_0>0$ is a generic positive constant. Using Young's inequality $a b\le \varepsilon a^2 + b^2/(4\varepsilon), \varepsilon>0,$ we can conclude that the linear operator $A_2$ (now considered as a linear operator from $E_1$ into $E_0$) has the relative zero norm, i.e. for any $\varepsilon>0$ there exists a constant $K_\varepsilon>0$ such that 
$\Vert A_2 \Phi\Vert_{E_0}\le \varepsilon \Vert\Phi\Vert_{E_1} + K_\varepsilon \Vert\Phi\Vert_{E_0}$.  (cf. \cite[Section 2.1]{Paz83} and also \cite{A4}). By virtue of \cite[Lemma 2.5]{A4} (see also \cite[Theorem 2.1]{Paz83}) the class ${\mathcal M}_1$ is closed with respect to perturbations by linear operators with zero relative norm. Thus ${\mathcal G}^\prime(\bar\Phi)\in {\mathcal M}_1(E_0,E_1)$ for any $\bar\Phi\in {\mathcal O}_1$. The proof of the short time existence of a solution $\Phi$  now follows from the aforementioned abstract result \cite[Theorem 2.7]{A4}.
\hfill$\diamondsuit$

\begin{remark}
\rm 
The nonlocal quantities $L = \int_\Gamma \ds,\ A = \frac12 \int_\Gamma \det(\vecx,\partial_s\vecx )\ds,\  {\mathcal E} = \int_\Gamma k^2\ds$ are $C^1$ mappings from a neighborhood ${\mathcal O}_{\frac{1}{2}}$ into $\R$. Indeed, as 
\[
L = \int_0^1 \vert\partial_u \vecx\vert \du,
\quad
A=\frac12 \int_0^1 \det(\vecx,\partial_u\vecx ) \du,
\quad
{\mathcal E} = \int_0^1 k^2 \vert\partial_u \vecx\vert \du,
\]
we have $L,A,{\mathcal E}\in  C^1({\mathcal O}_{\frac{1}{2}} ,\R )$. For the nonlocal term we conclude $\scF_\Gamma \in C^1({\mathcal O}_{\frac{1}{2}} ,\R )$ provided that 
$\scF_\Gamma = \scF (L,A,{\mathcal E})$ where $\scF (L,A,{\mathcal E})$ is a $C^1$ smooth function of its arguments $L,A, {\mathcal E}$. 
\end{remark}

\medskip
Let us introduce the $C^0$ and $C^1$ norms of a  quantity ${\sf F}:\G\to \R$ defined on a smooth  curve $\Gamma$ as follows
\begin{equation}
\Vert {\sf F} \Vert_{0,\Gamma}= \max_{\Gamma} |{\sf F}|, \qquad 
\Vert {\sf F} \Vert_{1,\Gamma}= \max_{\Gamma} \left( |{\sf F}| + |\partial_s {\sf F}| \right) .
\label{C1norm}
\end{equation}

\medskip
\begin{theorem}\label{globalexistence}
Let $\Gamma^t, t\ge 0,$ be a family of planar curves evolving in the normal direction with the velocity $\beta$ for which short time existence of smooth solutions is guaranteed by Theorem~{\rm\ref{localexistence}}. Suppose that $\Gamma^t, 0\le t<T_{max},$ is a maximal solution defined on the maximal time interval $[0, T_{max})$. If $T_{max}<+\infty$ then either $\Vert k\Vert_{0,\Gamma^t} + \Vert \beta \Vert_{0,\Gamma^t}\to \infty$ as $t\to T_{max}$ or $\sup_{0\le t<T_{max}} (\Vert k\Vert_{0,\Gamma^t} + \Vert \beta \Vert_{0,\Gamma^t}) < \infty$ and, in this case,  $\liminf_{t\to T_{max}^-} \min_{\Gamma^t} \tilde\beta^\prime_k(\vecx,k,\nu) = 0$.
\end{theorem}

\medskip
\noindent P r o o f:  Similarly as in \cite[Theorem 3.1]{A2} we can argue by a  contradiction. Suppose to the contrary that the sum $\Vert k\Vert_{0,\Gamma^t} + \Vert \beta \Vert_{0,\Gamma^t}$ of $C^0$ norms is bounded on $[0,T_{max})$ for $T_{max}<+\infty$ and  $\liminf_{t\to T_{max}^-} \min_{\Gamma^t}\beta^\prime_k > 0$. Following the continuation argument due to Gage \cite{Gage1986} and Angenent \cite{A2}, we shall prove that there exists a limiting curve $\Gamma^{T_{max}}=\Image(\vecx(.,T_{max}))$ where the $\vecx(u,T_{max}) = \lim_{t\to T_{max}} \vecx(u,t)$ uniformly w.r. to $u\in[0,1]$. Since the tangential velocity $\alpha$ has no impact on the shape of $\Gamma^t$ and the existence of the limiting curve $\G^{T_{max}}$ we can assume $\alpha=0$ in equation (\ref{curvatureeq}) for the curvature. With regard to the assumption on boundedness of the curvature $k$ and normal velocity $\beta$ we conclude that the reaction term $k^2\beta$ in the curvature equation remains uniformly bounded for $u\in[0,1]$ and $0\le t<T_{max}$. The parabolic equation (\ref{curvatureeq}) has therefore a smooth solution $k$ up to the limiting time $t=T_{max}$. Now, it follows from the position vector equation $\partial_t\vecx = \beta(\vecx, k,\nu) \vecN$ and the tangent angle equation (\ref{eq:equation-nu2}) that the limiting curve $\Gamma^{T_{max}}$ exists and is smooth. Moreover, using equation (\ref{eq:theta}) and expression $\theta=\ln\ratio_{\varphi}$ for the $\varphi$-adjusted local length $\ratio_{\varphi}$ we also conclude existence of the limit $\lim_{t\to T_{max}}\ratio_{\varphi}(u,t) = \ratio_{\varphi}(u,T_{max}) >0$ uniformly w.r. to $u\in[0,1]$. Starting from the limiting curve $\Gamma^{T_{max}}$ and using local in time existence result from Theorem~\ref{localexistence}, the family of curves $\Gamma^t$ can be prolonged to a larger time interval $[0,T')$ where $T'>T_{max}$. This is a contradiction. 
\hfill$\diamondsuit$

\begin{remark}
\rm 
In general, an estimate for $\Vert k\Vert_{0,\Gamma}$ is not sufficient in order to control behavior of the norm $\Vert \beta\Vert_{0,\Gamma}$ and vice versa. Indeed, let us consider, for example, $\beta = k - \vert \vecx \vert^2$. The selfsimilar family of circles parameterized by $\vecx(u,t) = R(t) (\cos2\pi u,\sin 2\pi u)^T$ evolves with the normal velocity $\beta$ provided that $-\dot{R}(t) = R(t)^{-1} - R(t)^2$. If $R(0)>1$ then $T_{max}<\infty$ and $\Vert k\Vert_{0,\Gamma^t} \to 0$ holds but  $\Vert \beta\Vert_{0,\Gamma^t} \to \infty$ as $t\to T_{max}$. On the other hand, if $\beta =\vert \vecx\vert^\delta k$, $1\le\delta <2$ then $\Vert \beta\Vert_{0,\Gamma^t}$ stays bounded but $\Vert k\Vert_{0,\Gamma^t} \to \infty$ holds as $t\to T_{max}$. Finally, if $\beta = (\vert \vecx\vert -1)^\delta k$, $0<\delta<1$, then   $-\dot{R}(t) = (R(t)-1)^{\delta}/R(t), R(0)>1,$ with $R(t)\to 1^+$ as  $t\to T_{max}<\infty$. Moreover, $\Vert k \Vert_{0,\Gamma^t}+\Vert \beta\Vert_{0,\Gamma^t}$ stays bounded for $t\in[0,T_{max})$ but $\liminf_{t\to T_{max}^-} \min_{\Gamma^t} \tilde\beta^\prime_k(\vecx,k,\nu) = 0$.
\end{remark}

\section{Affine scale space evolution and crystalline curvature tangential velocity}
\label{sec:affine}

In the theory of image processing, the so-called morphological image and shape multiscale analysis is often used because of its contrast and affine invariance properties. 
Affine scale space evolution of closed planar curves has been introduced and studied by Angenent, Sapiro and Tannenbaum \cite{AST1998,ST1994} and Alvarez \emph{et al.} \cite{AGL1993}. They derived a geometric equation  of the form (\ref{geomrov}) with the normal velocity given by
\begin{equation}
\beta(k)= k^{\frac13}.
\label{eq:affine}
\end{equation}
In this section we will show that the geometric equation (\ref{eq:affine}) representing affine scale space evolution is closely related to the crystalline curvature adjusted tangential velocity. For a wide class of tangential velocities, local existence, uniqueness and continuation of smooth solutions to the geometric equation (\ref{eq:affine}) has been shown by Mikula and \v{S}ev\v{c}ovi\v{c} in \cite{MikulaS2001}.

In order to investigate affine scale space evolution, we have to introduce a notion of the so-called affine arc-length parameterization (cf. \cite{ST1994}). Consider a planar closed convex curve $\Gamma$ parameterized by $\vecx$ with the arc-length parameterization $s$, i.e. $\Gamma=\{ \vecx(s), s\in [0,L]\}$. Such a curve can be re-parameterized by a new parameterization $S$ such that  $\Gamma=\{ \vecX(S), S\in [0,{\mathcal L}]\}$ where $\vecX(S) = \vecx(s)$ and $\dS = \vartheta \,\ds$. Here $\vartheta\ge 0$ is a nonnegative function defined on the curve $\Gamma$. 
The parameterization $S$ is called {\it affine arc-length} parameterization iff 
\[
|K| \equiv 1 \quad \hbox{on} \ \ \Gamma, \quad \hbox{where} \ \ 
K = \det(\partial_S \vecX, \partial^2_S \vecX), 
\]
(cf. \cite{AST1998}). In other words, $S$ is a parameterization for which the modulus of the affine curvature $K$ is constant along the curve $\Gamma$. It is well known (see \cite{ST1994}) that the parameterization $S$ is the affine arc-length iff $\vartheta = |k|^{\frac13}$. Indeed, since $\dS = \vartheta \,\ds$, we have 
\[
K = \det(\partial_S \vecX, \partial^2_S \vecX)
= \det(\vartheta^{-1}\partial_s \vecx, \vartheta^{-2}\partial^2_s \vecx - \partial_s\vartheta  \vartheta^{-3}\partial_s \vecx)
= \vartheta^{-3} \det(\partial_s \vecx, \partial^2_s \vecx) = \vartheta^{-3} k.
\]
Hence, $|K|\equiv 1$ on $\Gamma$ iff $\vartheta = |k|^{\frac13}$, i.e. $\dS = |k|^{\frac13} \,\ds$ on $\Gamma$, as claimed. 

Having a new parameterization $S$ of a family $\Gamma^t= \Image(\vecX(., t))$ of planar curves we can study evolution of $\Gamma^t, t\in [0,T],$ for which the affine intrinsic heat equation 
\begin{equation}
\partial_t\vecX = \partial^2_S\vecX
\label{eq:affineheat}
\end{equation}
is fulfilled. Since $\dS = \vartheta \,\ds$ we have
\[
\partial^2_S\vecX = 
\vartheta^{-1} \partial_s \left( 
\vartheta^{-1} \partial_s \vecx
\right) 
= \vartheta^{-2} \partial^2_s \vecx 
- \frac{\partial_s\vartheta}{\vartheta^3} \partial_s \vecx
= \frac{k}{\vartheta^2} \vecN  - \frac{\partial_s\vartheta}{\vartheta^3}  \vecT.
\]
Substituting $\vartheta=|k|^{\frac13}$ to the this equation we obtain the fact that the affine intrinsic heat equation  can be rewritten in the form of a geometric equation $\partial_t\vecX = \beta \vecN + \alpha \vecT$ with the normal velocity $\beta$ and tangential velocity $\alpha$ given by 
\[
\beta = k^{\frac13}, \qquad \alpha = - \frac{\partial_s\vartheta}{\vartheta^3} = 
- \frac13 \frac{\partial_s k}{k^{\frac53}}
=
- \frac{\partial_s\beta}{k}. 
\]
But this is just the tangential velocity $\alpha$ intrinsically built in the crystalline curvature flow  (cf. \cite{Yazaki20XXa}). It corresponds to the curvature adjusted tangential velocity with $\phi(k)\equiv |k|$.  This way we have shown the following proposition.

\begin{proposition}\label{affineproposition}
Let $S$ be the affine arc-length parameterization. Then the solution $\vecX=\vecX(S,t)$ to the affine intrinsic heat equation $\partial_t \vecX = \partial^2_S\vecX$ represents a flow of planar curves with the normal velocity $\beta=k^{\frac13}$ and the crystalline curvature adjusted tangential velocity $\alpha = -\frac{\partial_s\beta}{k}$. 
\end{proposition}

\begin{remark}
\rm
Recall that the evolution equation (\ref{eq:equation-nu}) for the tangent angle $\nu$ can be deduced from the fact $\partial_t\nu = \partial_s\beta +\alpha k$ (cf. \cite{MikulaS2001}). Therefore the curvature adjusted tangential velocity $\alpha$ given by $\alpha = - \partial_s\beta / k$ corresponds to $\partial_t\nu =0$ (see \cite{Yazaki20XXa}). Hence $\nu=\nu(u)$ is a function of $u$ variable only. In general, if a solution curve is strictly convex then the function $\nu=\nu(u)$ is invertible because $g^{-1} \partial_u\nu = \partial_s \nu = k >0$. Therefore we can reparameterize a convex curve by using the tangent angle $\nu$ as its new parameterization, i.e. $\partial_s = k \partial_\nu$. Then $\alpha = - \partial_s\beta / k = - \partial_\nu\beta$ and $\partial^2_s\beta = - k^2 \partial_\nu\alpha - \alpha\partial_s k$. Now, it follows from equation (\ref{eq:equation-k}) that 
\begin{equation}
\partial_t k = k^2( \partial^2_\nu \beta + \beta).
\label{eq:gaussparam}
\end{equation}
This is very useful evolution equation for strictly convex curves expressed in the so-called Gauss-parameterization  (see \cite{GH1986} in the case $\beta=k$).
\end{remark}

\section{Evolution of plane curves with a nonlocal normal velocity}
\label{sec:nonlocal}

In this section we discuss the role and importance of geometric equation (\ref{geomrovnonloc}) with a nonlocal term $\scF_\Gamma$ depending the entire curve $\Gamma$. The normal velocity is given by 
$\beta(\vecx, k, \nu) \equiv \tilde\beta(\vecx, k, \nu) + \scF_\Gamma$, where $\tilde\beta(\vecx, k, \nu) :\R^2\times\R\times \R\to \R$ is a local part of the normal velocity and $\scF_\Gamma = \scF(L, A, {\mathcal E})$ represents its nonlocal part depending on the global quantities $L, A, {\mathcal E}$ computed over the entire curve.

\subsection{Enclosed area- and total length-preserving nonlocal flows}

Recall the area evolution equation for a flow of embedded closed plane curves driven in normal direction with a normal velocity $\beta$:
\begin{equation}
\frac{d}{dt} A + \int_{\Gamma} \beta\, \,\ds = 0,
\label{eq:area}
\end{equation}
where $A=A^t$ is the area enclosed by the curve $\Gamma=\Gamma^t$ (see e.g. \cite{MikulaS2001}). 

Our first example is the enclosed area-preserving flow having the form of (\ref{geomrovnonloc}). Inserting $\beta(\vecx, k, \nu) \equiv \tilde\beta(\vecx, k, \nu) + \scF_\Gamma$ into (\ref{eq:area}) we obtain 
\[
0 = \frac{d}{dt} A + \int_{\Gamma} \beta \,\ds 
=  \frac{d}{dt} A + \int_{\Gamma} \tilde \beta \,\ds + \scF_\Gamma \int_{\Gamma} \,\ds  = \int_{\Gamma} \tilde \beta \,\ds + L \scF_\Gamma.
\]
Hence, setting  $\scF_\Gamma=- \<\tilde \beta\>=\frac{1}{L}\int_\Gamma \tilde\beta \,\ds$ where $L=\int_\Gamma \,\ds$, we conclude $\frac{d}{dt} A =0$ for the normal velocity given by 
$\beta= \tilde\beta - \<\tilde \beta\>$. Hence evolution according to such a normal velocity represents the area-preserving flow. The particular case $\tilde\beta\equiv k$ has been intensively studied e.g. by Gage \cite{Gage1986}. In this case we have
\[
\beta= k  - \< k \> = k - \frac{1}{L} \int_\Gamma k \,\ds = k - \frac{2\pi}{L}.
\]
Notice that the flow $\beta=k-\frac{2\pi}{L}$ is  the area-preserving and length decreasing in a suitable metric. Indeed, by the Cauchy-Schwartz inequality we obtain 
\begin{equation}
(2\pi)^2 =  \left(\int_\Gamma k \,\ds\right)^2 \le \int_\Gamma \,\ds \int_\Gamma k^2 \,\ds = L \int_\Gamma k^2 \,\ds
\label{CBS} 
\end{equation}
and so 
\[
\frac{d}{dt} L
=-\int_{\Gamma}k\beta\,\ds
=-\int_{\Gamma}k^2\,\ds
+\frac{2\pi}{L}\int_{\Gamma}k\,\ds
=-\int_{\Gamma}k^2\,\ds+\frac{(2\pi)^2}{L}
\leq 0. 
\]
Moreover, Gage proved that for the area-preserving flow we have $T_{max}=\infty$ and $\Gamma^t$ approaches a circle provided that the initial curve $\Gamma^0$ is convex (see \cite[Theorem 4.1]{Gage1986}). 

The total length-preserving flow represents a dual notion to the area-preserving flow. Our purpose is to construct a nonlocal term $\scF_\Gamma$ in such a way that $\partial_t L =0$. Using the total length equation (\ref{eq:length}) with the normal velocity having the form: $\beta=\tilde\beta + \scF_\Gamma$,  we obtain
\[
\frac{d}{dt} L = - \int_{\Gamma^t} k\beta \,\ds = 
- \int_{\Gamma} k\tilde\beta \,\ds - \scF_\Gamma \int_{\Gamma} k \,\ds
=
- \int_{\Gamma} k\tilde\beta \,\ds - 2\pi \scF_\Gamma.
\]
Therefore, for the normal velocity $\beta$ given by $\beta= \tilde\beta - \frac{L \< k \tilde\beta\>}{2\pi}$  we have $\frac{d}{dt} L =0$, i.e. the total length-preserving flow. In particular, by choosing the local part of the velocity $\tilde \beta \equiv  k$ we conclude
\[
\beta = k - \frac{{\mathcal E}}{2\pi},
\]
where ${\mathcal E} = \int_\Gamma k^2 \ds$ is the total elastic energy. In order to emphasize that the length-preserving flow is a dual notion to the area-preserving flow, we note that the flow of curves evolved by the velocity $\beta = k - \frac{{\mathcal E}}{2\pi}$ preserves the total length and enlarges the enclosed area. Indeed, using the inequality (\ref{CBS}) we obtain 
\[
\frac{d}{dt}A
= - \int_{\Gamma}\beta\,\ds 
= - \int_{\Gamma} \left( k - \frac{{\mathcal E}}{2\pi} \right) \,\ds 
= - 2\pi + \frac{L \int_\Gamma k^2\ds}{2\pi} \geq 0. 
\]
The total length preserving flow of convex curves was investigated recently by Mu and Zhu in \cite{Ma2012}. In their analysis they employed the Gauss parameterization (\ref{eq:gaussparam}) for evolution of convex curves in the plane. By using the arc-length parameterization we were able to extend the area increasing property for the case of Jordan curves in the plane.

\subsection{Gradient flow for the isoperimetric ratio}

Finally, we analyze a gradient flow for the isoperimetric ratio defined by:
\[
\Pi = \frac{L^2}{4\pi A}.
\]
Such a flow was first investigated by Jiang and Pan in \cite{Jiang2008}. Similarly as Mu and Zhu in \cite{Ma2012}, they employed the Gauss parameterization (\ref{eq:gaussparam}) for evolution of convex curves and so results of \cite{Jiang2008} can be applied to convex curves. In what follows we show that some of their results (e.g. area increasing property) can be directly generalized to the case of evolution of Jordan curves by using the arc-length parameterization. In our approach we employ the curvature equation (\ref{eq:equation-k}) in order to prove convexity preservation as well as area increasing (for Jordan curves) and length decreasing property (for convex curves). Moreover, in section~7 we present an example of of evolution of non-convex curves with temporal increase of the total length (see Fig.~\ref{fig:isogradient-gageinequality}).
Recall the well-known fact $\Pi\ge 1$ for any smooth closed embedded plane curve $\Gamma$ and the strict inequality $\Pi>1$ holds unless $\Gamma$ is a circle. Using relations (\ref{eq:length}) and (\ref{eq:area}) we obtain
\[
\frac{d}{dt}\Pi
= \frac{L \partial_t L}{2\pi A} - \frac{L^2 \partial_t A}{4\pi A^2}
= - \frac{L}{2\pi A} \int_{\Gamma}\left(k - \frac{L}{2 A}\right) \beta\,\ds.
\]
Hence, the flow driven in the normal direction with the normal velocity: 
\begin{equation}
\beta = k- \frac{L}{2A} 
\label{eq:isoperimnormal}
\end{equation}
represents a gradient flow for the isoperimetric ratio $\Pi$, 
i.e. $\partial_t \Pi <0$ for $\beta\not\equiv0$. 
Notice that $\beta$ is a normal velocity containing the nonlocal term $\scF_\Gamma=-L/(2A)$.
Furthemore, using the area equation (\ref{eq:area}) we obtain
\begin{equation}
\frac{d}{dt}A
= - \int_{\Gamma}\beta\,\ds 
= - \int_{\Gamma} k\,\ds + \frac{L}{2A} \int_{\Gamma} \,\ds 
= 2\pi (\Pi -1) \geq 0, 
\label{areaineq}
\end{equation}
because $\int_{\Gamma} k \,\ds = 2\pi, L= \int_{\Gamma} \,\ds$ and for the isoperimetric ratio we have $\Pi(\Gamma)\ge 1$. From inequality (\ref{areaineq}) we conclude that the enclosed area is nodecreasing with respect to time. 

Next, by following the convexity argument due to Gage, we shall prove that evolution of plane Jordan curves with the normal velocity given by (\ref{eq:isoperimnormal}) preserves convexity, i.e. $\Gamma^t$ is convex for  $t\ge t_0$ provided that $\Gamma^{t_0}$ is convex. A key point in the convexity preservation argument (cf.  Gage \cite{Gage1986}) is an application of a  parabolic comparison argument. Notice that a curve $\Gamma$ is convex if its curvature $k>0$. In the case $\beta = k + \scF_\Gamma$ equation (\ref{eq:equation-k}) for the curvature has the form:
\[
\partial_t k=\partial_s^2 k +\alpha\partial_s k+k^2 (k + \scF_\Gamma). 
\]
Let us denote by $K(t) = \min_{\Gamma^t} k(.,t)$, i.e. $K(t)$ is the minimal curvature over the curve $\Gamma^t$. Denote by $s^*(t)\in [0, L^t]$ the argument of minimum of $k$. It means that $K(t) = k(s^*(t), t)$. Then 
$\partial_s k(s^*(t), t)=0$ and $\partial_s^2 k(s^*(t), t)\ge 0$. Hence
\[
K'(t) \ge K^2(t) (K(t) +  \scF_\Gamma^t).
\]
Suppose that $K$ is a solution to this  ordinary differential inequality existing on some interval $[t_0, T)$ and such that $K(t_0)>0$. Then, it should be obvious that $K(t)>0$ for $t\in[t_0, T)$ provided that 
\begin{equation}
\inf_{t_0\le t < T}\scF_\Gamma^t > -\infty.
\label{Fbound}
\end{equation}
In order to prove convexity preservation for the isoperimetric ratio  gradient flow it is therefore sufficient to verify that the nonlocal part $\scF_\Gamma^t=-L^t/(2A^t)$ remains bounded from below for $t\ge t_0$. Indeed, as $\beta = k -L/(2A)$ we obtain 
\[
\frac{d}{dt} \frac{L}{A} = \frac{\partial_t L}{A} - \frac{L\partial_t A}{A^2} = -\frac{1}{A} \int_\Gamma k\beta\, \,\ds  + \frac{L}{A^2} \int_\Gamma \beta\, \,\ds
= -\frac{1}{A} \int_\Gamma \beta^2\, \,\ds  + \frac{L}{2 A^2} \int_\Gamma \beta\, \,\ds \le 0,
\]
since $\int_\Gamma \beta\, \,\ds = -\partial_t A\le 0 $ by (\ref{areaineq}). Therefore, 
\begin{equation}
- \frac{L^{t_0}}{2 A^{t_0}} < -\frac{L^t}{2 A^t} \le 0 \quad\text{for any}\  t\ge t_0, 
\label{ineqLA}
\end{equation}
and so condition (\ref{Fbound}) is satisfied by the the isoperimetric ratio gradient flow (\ref{eq:isoperimnormal}). 

Another important property of this flow can be deduced from the well known isoperimetric inequality 
\begin{equation}
\pi \frac{L}{A} \le \int_\Gamma k^2 \,\ds
\label{gageineq}
\end{equation}
derived by Gage \cite{Gage1983} for any convex $C^2$ smooth embedded closed curve in $\R^2$. Therefore, if $\Gamma^{t_0}$ is convex at some time $t_0\ge 0$, then, with regard to (\ref{eq:length}), we obtain 
\[
\frac{d}{dt} L^t = -  \int_{\Gamma^t} k\beta \,\ds  = -  \int_{\Gamma^t} k^2 \,\ds
+ \frac{L^t}{2A^t}  \int_{\Gamma^t} k \,\ds  = -  \int_{\Gamma^t} k^2 \,\ds
+ \pi \frac{L^t}{A^t}  \le 0 \ \ \text{at}\ t=t_0.
\]
It means that the isoperimetric ratio gradient flow  is a curve-shortening flow of evolving convex curves.

Summarizing result from this subsection we can conclude the following proposition.

\begin{proposition}\label{isoproposition}
A flow of planar embedded closed curves driven in the normal direction by the velocity $\beta = k- \frac{L}{2A}$ is a gradient flow minimizing the isoperimetric ratio $\Pi = L^2/(4\pi A)$. It is a convexity-preserving flow. It is the enclosed area increasing flow. Finally, it represents the total length decreasing flow for convex curves. 
\end{proposition}

Note that the convexity preservation and length decreasing property have been already established by Jiang and Pan in \cite{Jiang2008} for evolution of convex curves. Since we have adopted the arc-length paramterization, the area increasing property holds true for evolution of Jordan curves evolved in the normal direction by the velocity (\ref{eq:isoperimnormal}).

\section{Numerical scheme}
\label{sec:scheme}



The purpose of this section is to construct a numerical approximation scheme for 
solving the system of governing equations (\ref{eq:equation-k})--(\ref{eq:equation-x}) complemented by  the curvature adjusted tangential velocity equation (\ref{eq:tangential_velocity}) satisfying the renormalization constraint $\<\varphi \alpha\>=0$.

{\bf Discretization setting.}\ 
Given an initial $N$-sided polygonal curve $\calP^0=\bigcup_{i=1}^N\calS_i^0$, 
our purpose is to find a family of $N$-sided polygonal curves $\{\calP^j\}_{j=1, 2, \cdots}$, 
$\calP^j=\bigcup_{i=1}^N\calS_i^j$, where $\calS_i^j=[\vecx_{i-1}^j, \vecx_{i}^j]$ is the $i$-th edge with $\vecx_{0}^j=\vecx_N^j$ for $j=0, 1, 2, \cdots$. 
The initial polygon $\calP^0$ is an approximation of $\Gamma^0$ satisfying 
$\{\vecx_i^0\}_{i=1}^N\subset\calP^0\cap\Gamma^0$, 
and $\calP^j$ is an approximation of $\Gamma^t$ at the time $t=t_j$, 
where $t_j=j\tau$ is the $j$-th discrete time level, $j=0, 1, 2, \cdots$,  if we use a fixed time increment $\tau>0$, or $t_j=\sum_{l=0}^{j-1}\tau_l$ if we use 
adaptive time increments $\tau_l>0, l=0, \cdots, j-1 $. 
The updated curve $\calP^{j+1}$ is determined from the data for $\calP^j$ at the previous time step by 
using the following numerical discretization in space and time:  

{\bf Numerical algorithm.}\ 
A flowchart of the numerical algorithm can be stated as follows. 
Hereafter, we will use ``${\sf func}$'' for a given quantity or function. 
Note that local quantities $\{\vecx_i\}$, $\{k_i\}$, $\{p_i\}$, $\{\alpha_i\}$ should  satisfy the periodic conditions: 
${\sf func}_0={\sf func}_N$, 
${\sf func}_{N+1}={\sf func}_{1}$. As for the tangent angle vector $\{\nu_i\}$ we have mod-$2\pi$ boundary conditions: $\nu_0\equiv\nu_N$ (mod-$2\pi$),  $\nu_{N+1}\equiv \nu_{1}$ (mod-$2\pi$).

\begin{description}
\item[\it Step 0.]\ 
Put the initial step $j=0$; 
\item[\it Step 1.]\ 
Construct the data at the  $j$-th time step as follows: 
\begin{description}
\item[\it Step 1-1.]\ 
Using $\{\vecx_i^j\}$, compute the local quantities $\{p_i^j\}$, $\{\nu_i^j\}$, $\{k_i^j\}$, 
nonlocal quantities $L^j$, $A^j$, ${\mathcal E}^j$ 
and $\{\beta_i^j=\tilde{\beta}(\vecx_i^j, \nu_i^j, k_i^j)+\scF_{\calP^j}\}$; 
\item[\it Step 1-2.]\ 
Solve equation  (\ref{eq:tangential_velocity}) for $\alpha$ and compute  $\{\alpha_i^j\}$; 
\end{description}
\item[\it Step 2.]\ 
Compute updated data $\{\bar{p}_i\}$, $\{\bar{\nu}_i\}$, $\{\bar{k}_i\}$, $\bar{L}$, as follows: 
\begin{description}
\item[\it Step 2-1.]\ 
Solve equation (\ref{eq:equation-g}) for $g$ and compute  $\{\bar{p}_i\} = N g_i^j$, $\bar{L} =\sum_{i=1}^N \bar{p}_i$; 
\item[\it Step 2-2.]\ 
Solve equation (\ref{eq:equation-k}) for $k$ and compute  $\{\bar{k}_i\}$, $\<{\sf func}(\bar{k})\>$; 
\item[\it Step 2-3.]\ 
Solve equation (\ref{eq:equation-nu}) for $\nu$ and compute  $\{\bar{\nu}_i\}$; 
\end{description}
\item[\it Step 3.]\ 
Solve equation (\ref{eq:equation-x}) for the position vector $\vecx$ and compute  $\{\vecx_i^{j+1}\}$; 
\item[\it Final step.]\ 
Put $j:=j+1$, go to {\it Step 1}, and compute  all updated quantities. 
\end{description}

In what follows, we shall explain the aforementioned steps in more details. 
In the following explanation of Step 1-1 and Step 1-2, 
we will omit the superscript $j$. 

{\bf Step 1-1}.\ 
From the $j$-th step data $\{\vecx_i^j\}$, 
we will define local quantities $\{p_i^j\}$, $\{\nu_i^j\}$, $\{k_i^j\}$, 
and nonlocal quantities $L^j$, $A^j$, $\<{\sf func}(k^j)\>, \cdots$. 

{\bf Construct $\bm{\{p_i^j\}}$ and $\bm{L^j}$}.\ 
The length of $\calS_i$ is given by $p_i=|\vecx_{i}-\vecx_{i-1}|$. 
Then we have the total length of $\calP$ such as 
$L=\sum_{i=1}^Np_i$. 

{\bf Construct $\bm{\{\nu_i^j\}}$ and $\bm{A^j}$}.\ 
The $i$-th unit tangent vector $\vecT_i$ is defined by 
$\vecT_i=\frac{1}{p_i}(\vecx_{i}-\vecx_{i-1})$. 
Then the $i$-th unit tangent angle $\nu_i$ is obtained from 
$\vecT_i=(\cos\nu_i, \sin\nu_i)^{\mathrm{T}}$ in the following way: 
firstly, from $\vecT_1=(T_{11}, T_{21})^{\mathrm{T}}$, 
we obtain $\nu_1=2\pi-\arccos(T_{11})$ if $T_{12}<0$; 
$\nu_1=\arccos(T_{11})$ if $T_{12}\geq 0$. 
Secondly, for $i=1, 2, \cdots, N$,  we successively compute $\nu_{i+1}$ from $\nu_{i}$ as follows: 
\[
\nu_{i+1}=\left\{\begin{array}{@{}ll}
\nu_i+\arcsin(D), & \mbox{if $I>0$}, \\
\nu_i+\arccos(I), & \mbox{if $D>0$}, \\
\nu_i-\arccos(I), & \mbox{otherwise}, 
\end{array}\right.
\ 
\mbox{where}\ D=\det(\vecT_i, \vecT_{i+1}), 
\ 
I=\vecT_i\cdot\vecT_{i+1}.
\]
Finally, we obtain 
$\nu_0=\nu_1-(\nu_{N+1}-\nu_{N})$ and 
$\nu_{N+2}=\nu_{N+1}+(\nu_2-\nu_1)$. 
If needed, we can also calculate the enclosed area by formula $A=-\frac{1}{2}\sum_{i=0}^N\vecx_i\cdot\vecN(\nu_i)p_i$ with the inward normal vector 
$\vecN_i=(-\sin\nu_i, \cos\nu_i)^{\mathrm{T}}$.

{\bf Construct $\bm{\{k_i^j\}}$ and $\bm{\<{\sf func}(k^j)\>}$}.\ 
In order to derive a discrete numerical scheme, 
we follow the flowing finite volume method adopted for curve evolutionary problems as it was proposed by 
Mikula and \v{S}ev\v{c}ovi\v{c} in~\cite{MikulaS2004b}. 
Let us introduce the dual volume 
$\calS_i^*=[\vecx_{i}^*, \vecx_{i}]\cup[\vecx_{i}, \vecx_{i+1}^*]$ of $\calS_i$, 
where $\vecx_{i}^*=\frac{1}{2}(\vecx_{i}+\vecx_{i-1})$ for $i=1, 2, \cdots, N$. 
We define the $i$-th unit tangent angle of $\calS_i^*$ by $\nu_i^*=\frac{1}{2}(\nu_i+\nu_{i+1})$. 
The $i$-th curvature $k_i$ is approximated by a constant value on $\calS_i$, which 
is obtained from integration of $k=\partial_s\nu$ over $\calS_{i}$ with respect to $s$: 
\[
\int_{\calS_{i}}k\,ds
=k_i\int_{\calS_{i}}\,ds
=k_ip_i, \quad
\int_{\calS_{i}}k\,ds
=\int_{\calS_{i}}\partial_s\nu\,ds
=[\nu]_{\vecx_{i-1}}^{\vecx_{i}}
=\nu_i^*-\nu_{i-1}^*. 
\]
Hereafter, 
$\int_{\calS_{i}}{\sf func}\,\ds$ means 
$\int_{s_{i-1}}^{s_{i}}{\sf func}\,\ds$ for the arc-length $s_i$ parameter satisfying $\vecx_{i}=\vecx(s_i, \cdot)$. 
Thus we have $k_i=(\partial_s\nu^{*})_i$, where 
$(\partial_s{\sf func})_i=({\sf func}_{i}-{\sf func}_{i-1})/p_{i}$. 
The $i$-th curvature $k_i^*$ at $\vecx_i$ can be defined as 
$k_i^*=\frac{1}{2}(k_{i}+k_{i+1})$. The curvature is assumed to be  constant $k_i^*$ over the finite volume  $\calS_i^*$. 
We also compute $\<{\sf func}(k)\>=\frac{1}{L}\sum_{i=1}^N{\sf func}(k_i)p_i$. 

{\bf Step 1-2 (construct  $\bm{\{\alpha_i^j\}}$)}.\ 
We discretize equation (\ref{eq:tangential_velocity}) for the tangential velocity $\alpha$: 
\[
\partial_s(\varphi\alpha)
=\frac{\<f\>}{\<\varphi\>}\varphi-f+\left(\frac{L}{g}\<\varphi\>-\varphi\right)\omega, 
\]
where $f=(\partial_s^2\beta+k^2\beta)\varphi'(k)-k\beta\varphi(k)$. 
Integrating the previous equation over $\calS_i$ yields: 
\begin{align*}
\int_{\calS_i}\partial_s(\varphi\alpha)\,ds
&=\big[\varphi(k)\alpha\big]_{\vecx_{i-1}}^{\vecx_{i}} \\\
&=\frac{\<f\>}{\<\varphi\>}\int_{\calS_i}\varphi(k)\,ds
-\int_{\calS_i}f\,ds
+\left(L\<\varphi\>\int_{\calS_i}\frac{1}{g}\,ds-\int_{\calS_i}\varphi(k)\,ds\right)\omega.
\end{align*}
Hence
\begin{align*}
\psi_i
&=\varphi(k_i^*)\alpha_i-\varphi(k_{i-1}^*)\alpha_{i-1} =\frac{\<f\>}{\<\varphi\>}\varphi(k_i)p_i
-f_ip_i+\left(L\<\varphi\>\frac{1}{N}-\varphi(k_i)p_i\right)\omega, \\
f_i
&=\big((\partial_s(\partial_{s^*}\beta))_i
+k_i^2\beta_i\big)\varphi'(k_i)
-k_i\beta_i\varphi(k_i), 
\end{align*}
where $\beta_i=\tilde{\beta}(\vecx_i^*, \nu_i, k_i)+\scF_{\calP}$ is assumed to be constant on the finite volume $\calS_i$, 
\[
(\partial_s(\partial_{s^*}\beta))_i
=\frac{(\partial_{s^*}\beta)_i-(\partial_{s^*}\beta)_{i-1}}{p_i}
=\frac{1}{p_i}\big[\partial_s\beta\big]_{\vecx_{i-1}}^{\vecx_{i}}, 
\quad
(\partial_{s*}{\sf func})_i=\frac{{\sf func}_{i+1}-{\sf func}_{i}}{p_{i}^*}, 
\]
$p_i^*=\frac{1}{2}(p_i+p_{i+1})$ is the length of $\calS_i^*$. 
The curve averages $\<f\>$ and $\<\varphi(k)\>$ are approximated as follows: 
\[
\<f\>
=\frac{1}{L}\sum_{i=1}^Nf_ip_i, \quad
\<\varphi(k)\>
=\frac{1}{L}\sum_{i=1}^N\varphi(k_i)p_i. 
\]
In order to uniquely determine $\{\alpha_i\}_{i=1}^N$ we have to take into account the renormalization constraint $\<\varphi(k)\alpha\>=0$. The average $\<\varphi(k)\alpha\>$ can be discretized as follows: 
\[
\<\varphi\alpha\>
\approx \frac{1}{L}\sum_{i=1}^N\varphi(k_i^*)\alpha_ip_i^*. 
\]
Notice that $L=\sum_{i=1}^Np_i=\sum_{i=1}^Np_i^*$. Let us define a partial sum of $\{\psi_i\}$ by
\[
\Psi_i=\sum_{l=2}^i\psi_l,\quad
i=2, 3, \cdots, N, \quad
\Psi_1=0. 
\]
Using this notation we obtain $\varphi(k_i^*)\alpha_ip_i^* =\varphi(k_1^*)\alpha_1p_i^*+\Psi_ip_i^*$. By summing these terms for $i=1, 2, \cdots, N$ we obtain
\[
\sum_{i=1}^N\varphi(k_i^*)\alpha_ip_i^*
=\varphi(k_1^*)\alpha_1 L
+\sum_{i=1}^N\Psi_ip_i^*, \quad
L=\sum_{i=1}^Np_i^*. 
\]
This way we have computed the vector $\{\alpha_i\}$ of tangential velocities:
\[
\alpha_1
=-\frac{1}{L\varphi(k_1^*)}\sum_{i=2}^N\Psi_ip_i^*, \quad
\alpha_i
=\frac{1}{\varphi(k_i^*)}
\left(\varphi(k_1^*)\alpha_1+\Psi_i\right), \quad
i=2, 3, \cdots, N. 
\]

In the next Step 2, 
the $j$-th step data $\{p_i^j\}$, $\{\nu_i^j\}$, $\{k_i^j\}$, $L^j$, $\<{\sf func}(k^j)\>$ are updated 
to the new data vectors $\{\bar{p}_i\}$, $\{\bar{\nu}_i\}$, $\{\bar{k}_i\}$, $\bar{L}$, $\<{\sf func}(\bar{k})\>$. Notice that these data vectors form  intermediate step quantities between $j$ and $j+1$-th time level.  
In order to update the position vector in Step 3, and they are updated to the $(j+1)$-th step data by means of $\{\vecx_i^{j+1}\}$ as in the final step. \footnote{If one wants to construct a simpler and faster (but less accurate) scheme Step 2 can be skipped by following discretization scheme proposed by the authors in \cite{SY2008}. }

{\bf Step 2-1 (compute $\bm{\{\bar{p}_i\}}$ and $\bm{\bar{L}}$)}.\ 
We semi-implicitly discretize equation (\ref{eq:equation-g}) for $g$ over $\calS_i^j$ to obtain
\[
\frac{\bar{p}_i-p_i^j}{\tau}
=-k_i^j\beta_i^j\bar{p}_i+\alpha_i^j-\alpha_{i-1}^j, \quad
i=1, 2, \cdots, N. 
\]
It can be solved as follows: 
\[
\bar{p}_i
=\frac{1+(\partial_s\alpha)_i^j\tau}{1+k_i^j\beta_i^j\tau}p_i^j
=\frac{p_i^j+(\alpha_i^j-\alpha_{i-1}^j)\tau}{1+k_i^j\beta_i^j\tau}, \quad
i=1, 2, \cdots, N. 
\]
Note that $g=|\partial_u\vecx|$ on $\calS_i^j$ is discretized as 
$Np_i^j=\frac{|\vecx_i^j-\vecx_{i-1}^j|}{1/N}$. 
Then we obtain $\bar{L}=\sum_{i=1}^N\bar{p}_i$. 

{\bf Step 2-2 (compute $\bm{\{\bar{k}_i\}}$ and $\bm{\<{\sf func}(\bar{k})\>}$)}.\ 
We semi-implicitly discretize equation (\ref{eq:equation-k}) for $k$ on $\calS_i^j$. For $i=1, 2, \cdots, N$, we obtain
\begin{align*}
\frac{\bar{k}_i-k_i^j}{\tau}
=\partial_{\bar{s}}((\tilde{\beta}_k')^{j}(\partial_{\bar{s}^*}\bar{k})
+(\tilde{\beta}_{\nu}')^{j}\bar{k}^*
+(\nabla_{\vecx}\tilde{\beta})^{j}\cdot\vecT(\nu^{*j}))_i
+\alpha_i^{*j}(\partial_{\bar{s}}\bar{k}^*)_i
+(k_i^j)^2\beta_i^{\bar{j}}.
\end{align*}
We and up with the tridiagonal system: 
\[
a_{-}\bar{k}_{i-1}+a_{0}\bar{k}_{i}+a_{+}\bar{k}_{i+1}
=k_i^j+(\partial_{\bar{s}}((\nabla_{\vecx}\tilde{\beta})^{j}\cdot\vecT(\nu^{*j}))_i
+(k_i^j)^2\beta_i^{\bar{j}})\tau, \quad
i=1, 2, \cdots, N, 
\]
subject to periodic boundary conditions:  
$\bar{k}_0=\bar{k}_N$, $\bar{k}_{N+1}=\bar{k}_1$
where
\[
a_{-}=-\frac{\tau}{\bar{p}_i}\left(
\frac{(\tilde{\beta}_{k}')_{i-1}^{j}}{\bar{p}_{i-1}^*}
-\frac{(\tilde{\beta}_{\nu}')_{i-1}^{j}}{2}
-\frac{\alpha_i^{*j}}{2}\right), \quad
a_{+}=-\frac{\tau}{\bar{p}_i}\left(
\frac{(\tilde{\beta}_{k}')_{i}^{j}}{\bar{p}_i^*}
+\frac{(\tilde{\beta}_{\nu}')_{i}^{j}}{2}
+\frac{\alpha_i^{*j}}{2}\right), 
\]
and $a_0=1-(a_{-}+a_{+})$. 
Here the values of 
$(\tilde{\beta}_{k}')_{i}^{j}$, 
$(\tilde{\beta}_{\nu}')_{i}^{j}$, 
$(\nabla_{\vecx}\tilde{\beta})_{i}^{j}\cdot\vecT(\nu_i^{*j})$ 
are evaluated at $(\vecx_i^{*j}, \nu_i^j, k_i^j)$ and 
$\beta_i^{\bar{j}}=\tilde{\beta}(\vecx_i^{*j}, \nu_i^j, k_i^j)+\scF(\bar{L}, A^j, \<{\sf func}(k^{\bar{j}})\>, \cdots)$, 
where $\<{\sf func}(k^{\bar{j}})\>=\frac{1}{\bar{L}}\sum_{i=1}^N{\sf func}(k_i^j)\bar{p}_i$. 
We also obtain updated nonlocal data 
$\<{\sf func}(\bar{k})\>=\frac{1}{\bar{L}}\sum_{i=1}^N{\sf func}(\bar{k}_i)\bar{p}_i$. 

{\bf Step 2-3 (compute $\bm{\{\bar{\nu}_i\}}$)}.\ 
We semi-implicitly discretize equation (\ref{eq:equation-k}) for $\nu$ over $\calS_i^j$. We obtain
\[
\frac{\bar{\nu}_i-\nu_i^j}{\tau}
=(\tilde{\beta}_k')_i^{\bar{j}}(\partial_{\bar{s}}(\partial_{\bar{s}^*}\bar{\nu}))_i
+(\alpha_i^{*j}+(\tilde{\beta}_\nu')_i^{\bar{j}})(\partial_{\bar{s}}\bar{\nu}^*)_i
+(\nabla_{\vecx}\tilde{\beta})_i^{\bar{j}}\cdot\vecT(\nu_i^j), \quad
i=1, 2, \cdots, N, 
\]
where the boundary data are given by 
$\bar{\nu}_{0}=\bar{\nu}_{2}-2\bar{k}_1\bar{p}_1$, 
$\bar{\nu}_{N+1}=\bar{\nu}_{N-1}+2\bar{k}_{N}\bar{p}_{N}$, 
which are the same as $\bar{k}_i=(\partial_{\bar{s}}\bar{\nu}^*)_i$ for $i=1, N$, respectively. 
It follows from the above discretized equation that $\nu$ can be computed by solving the tridiagonal system: 
\[
a_{-}\bar{\nu}_{i-1}+a_{0}\bar{\nu}_{i}+a_{+}\bar{\nu}_{i+1}
=b_i, \quad
b_i=\nu_i^j+(\nabla_{\vecx}\tilde{\beta})_i^{\bar{j}}\cdot\vecT(\nu_i^{j})\tau, \quad
i=1, 2, \cdots, N, 
\]
since the first and the $N$-th rows are given by
\[
a_{0}\bar{\nu}_{1}+(a_{-}+a_{+})\bar{\nu}_{2}=b_1+2a_{-}\bar{k}_1\bar{p}_1, \quad
(a_{-}+a_{+})\bar{\nu}_{N-1}+a_{0}\bar{\nu}_{N}=b_N-2a_{+}\bar{k}_N\bar{p}_N, 
\]
respectively, 
where
\[
a_{-}=-\frac{\tau}{\bar{p}_i}\left(
\frac{(\tilde{\beta}_k')_i^{\bar{j}}}{\bar{p}_{i-1}^{*}}
-\frac{\alpha_i^{*j}+(\tilde{\beta}_{\nu}')_i^{\bar{j}}}{2}
\right), \quad
a_{+}=-\frac{\tau}{\bar{p}_i}\left(
\frac{(\tilde{\beta}_k')_i^{\bar{j}}}{\bar{p}_{i}^{*}}
+\frac{\alpha_i^{*j}+(\tilde{\beta}_{\nu}')_i^{\bar{j}}}{2}
\right), 
\]
and $a_0=1-(a_{-}+a_{+})$. 
Here the values of 
$(\tilde{\beta}_k')_i^{\bar{j}}$, 
$(\tilde{\beta}_{\nu}')_i^{\bar{j}}$ and 
$(\nabla_{\vecx}\tilde{\beta})_i^{\bar{j}}$ are evaluated at 
$(\vecx_i^{*j}, \nu_i^j, \bar{k}_i)$, respectively. 
Having constructed the new tangential angle vector $\{\bar{\nu}_i\}_{i=1}^N$, 
we extend it for $i=0,1,N,N+1$ as follows: 
$\bar{\nu}_{0}=\bar{\nu}_{2}-2\bar{k}_1\bar{p}_1$, 
$\bar{\nu}_{N+1}=\bar{\nu}_1+(\bar{\nu}_{N}-\bar{\nu}_{0})$, 
and also 
$\bar{\nu}_{N+2}=\bar{\nu}_{N+1}+(\bar{\nu}_2-\bar{\nu}_1)$ if it's necessary. 

{\bf Step 3 (compute $\bm{\{\vecx_i^{j+1}\}}$)}.\ 
Now we are ready to obtain updated values of $\{\vecx_i^{j+1}\}$ by solving semi-implicitly discretized 
equation (\ref{eq:equation-x}) for the position vector $\vecx$ over  the finite volume $\calS_i^{*j}$:
\[
\frac{\vecx_i^{j+1}-\vecx_i^j}{\tau}
=w_i^{*\bar{j}}(\partial_{\bar{s}^*}(\partial_{\bar{s}}\vecx))_i^{j+1}
+\alpha_i^{j}(\partial_{\bar{s}^*}\vecx^*)_i^{j+1}
+(F_i^{*\bar{j}}+\scF_i^{\bar{j}})\vecN(\bar{\nu}_i^*), \quad
i=1, 2, \cdots, N. 
\]
We obtain the tridiagonal system with the periodic boundary condition 
$\bar{\vecx}_0=\bar{\vecx}_N$, $\bar{\vecx}_{N+1}=\bar{\vecx}_1$: 
\[
a_{-}\bar{\vecx}_{i-1}+a_{0}\bar{\vecx}_{i}+a_{+}\bar{\vecx}_{i+1}
=\vecx_i^j+(F_i^{*\bar{j}}+\scF_i^{\bar{j}})\vecN(\bar{\nu}_i^*)\tau, \quad
i=1, 2, \cdots, N, 
\]
where
\[
a_{-}=-\frac{\tau}{\bar{p}_i^*}\left(
\frac{w_i^{*\bar{j}}}{\bar{p}_{i}}
-\frac{\alpha_i^{j}}{2}\right), \quad
a_{+}=-\frac{\tau}{\bar{p}_i^*}\left(
\frac{w_i^{*\bar{j}}}{\bar{p}_{i+1}}
+\frac{\alpha_i^{j}}{2}\right), 
\]
and $a_0=1-(a_{-}+a_{+})$. 
Here we have used the following values: 
$w_i^{*\bar{j}}=w(\vecx_i^j, \bar{\nu}_i^*, \bar{k}_i^*)$, 
$F_i^{*\bar{j}}=F(\vecx_i^j, \bar{\nu}_i^*)$ and 
$\scF_i^{\bar{j}}=\scF(\bar{L}, A^j, {\mathcal E}^j)\>, \cdots)$.



\section{Computational results}
\label{sec:results}


In this section we present several computational examples of applications of the numerical scheme proposed in the previous section. In all our experiments we chose $N=100$ spatial grid points in order to discretize evolving curves. In all examples the initial curve has large variations in the curvature. Therefore accurate resolution of the evolved curve near narrow tails and kinks was required. Especially, in the case of thin dumb-bell initial curve connecting two open circular bubbles by a thin and long tubular channel, a very fine resolution of the curve and its curvature near the smooth connection of  circles and tube was required. We turn on the tangential redistribution by choosing the curvature adjusted tangential velocity $\alpha$ with the shape function $\varphi(k)$ of the form
\[
\varphi(k)=1-\eps+\eps\sqrt{1-\eps+\eps k^2} \quad \hbox{where}\  \eps=0.1. 
\]

\begin{figure}
\begin{center}
\includegraphics[width=0.3\textwidth,bb=120 74 334 278]{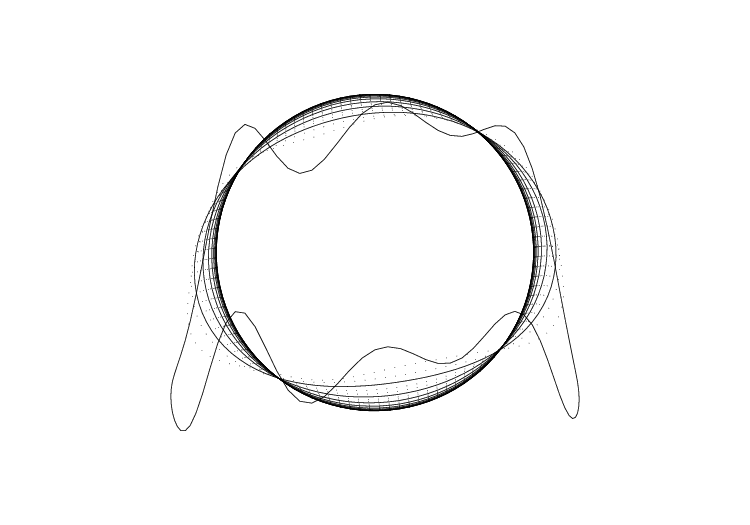}
\includegraphics[width=0.3\textwidth,bb=147 74 314 278]{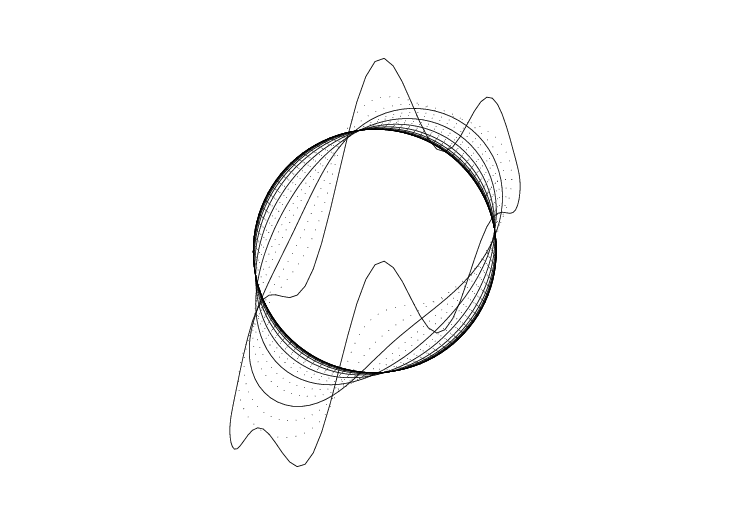}
\includegraphics[width=0.3\textwidth,bb=115 74 344 278]{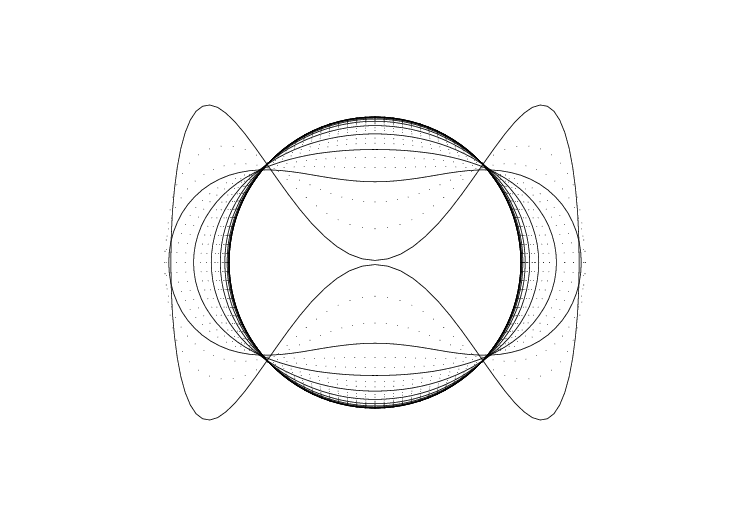}
\end{center}

\caption{Evolution of three initial curves by the area-preserving flow.}
\label{fig:areapreserving}
\end{figure}

In our first three examples shown in Fig.~\ref{fig:areapreserving} we present the area-preserving curve evolution with the normal velocity 
 \[
\beta= k - \frac{2\pi}{L}.
\]
We can observe stable computation of evolving curves approaching a circle with the area identical to the area of initial curve. 

\begin{figure}
\begin{center}
\includegraphics[width=0.52\textwidth,bb=126 160 334 196]{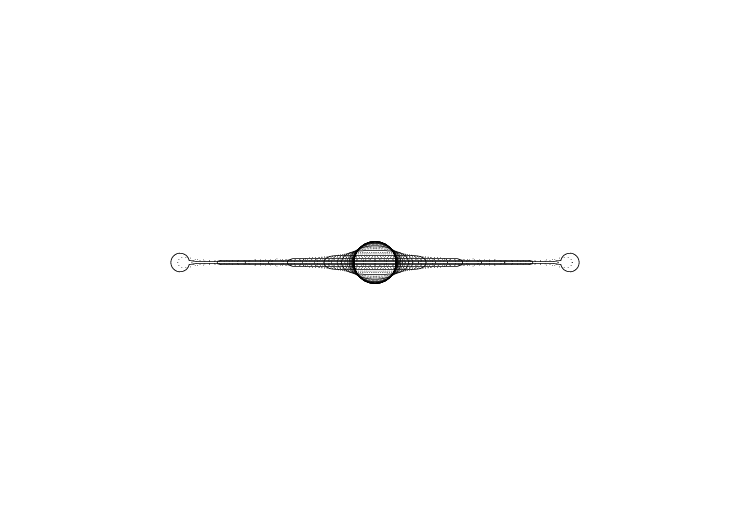}
\end{center}

\caption{Evolution of a thin dumb-bell initial curve by the area-preserving flow.}
\label{fig:areapreserving-2}
\end{figure}

In the next example depicted in Fig.~\ref{fig:areapreserving-2} we show evolution of a thin dumb-bell initial curve. Due to application of the curvature adjusted tangential velocity we obtained fine and accurate resolution of parts of the evolved curve having large modulus of the curvature and we could compute its evolution preserving the enclosed area over sufficiently large time interval until it became convex. The limiting curve is again a circle.

\begin{figure}
\begin{center}
\includegraphics[width=0.3\textwidth,bb=120 60 347 278]{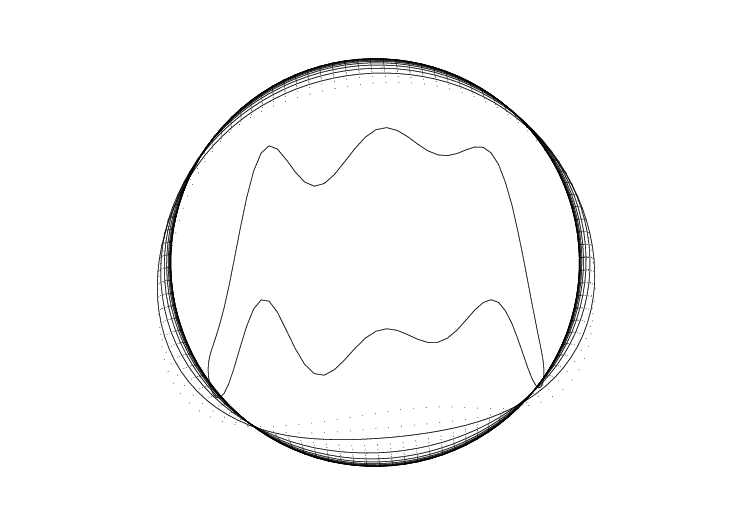}
\includegraphics[width=0.3\textwidth,bb=120 60 347 278]{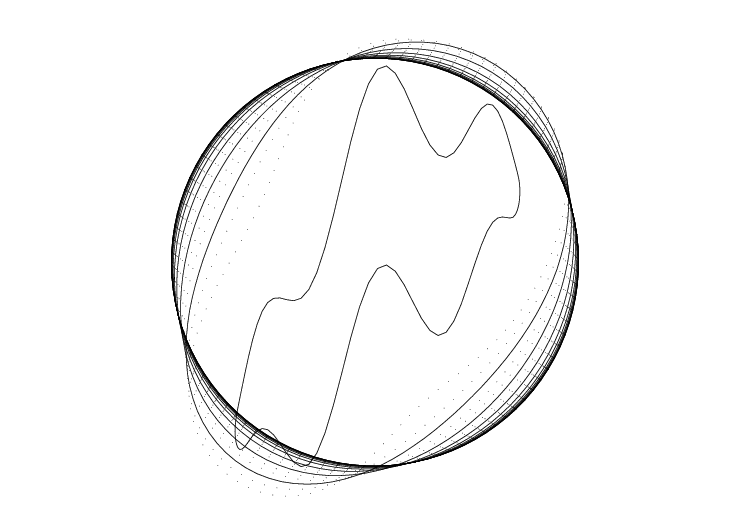}
\includegraphics[width=0.3\textwidth,bb=120 60 347 278]{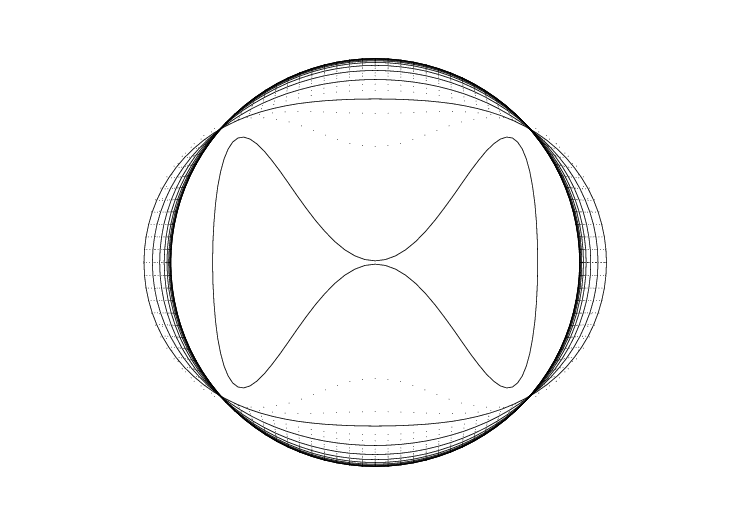}
\end{center}

\caption{Evolution of three initial curves by the total length-preserving flow. }
\label{fig:lengthpreserving}
\end{figure}

\begin{figure}
\begin{center}
\includegraphics[width=0.52\textwidth]{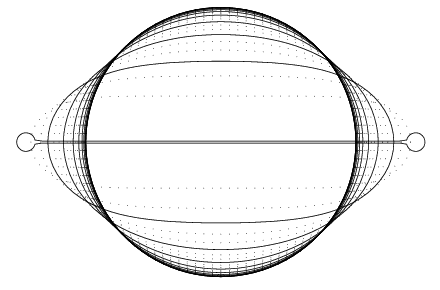}
\end{center}

\caption{Evolution of a thin dumb-bell initial curve by the total length-preserving flow. }
\label{fig:lengthpreserving-2}
\end{figure}

The second set of examples is shown in Fig.~\ref{fig:lengthpreserving}. It is  devoted to the total length-preserving flow with the normal velocity given by 
 \[
\beta= k - \frac{\mathcal E}{2\pi},
\]
where ${\mathcal E}=\int_\Gamma k^2 \ds$ is the total elastic energy. 
In Fig.~\ref{fig:lengthpreserving-2} evolution of the  thin dumb-bell initial curve is shown. In contrast to the area-preserving flow of the thin dumb-bell initial curve (see  Fig.~\ref{fig:areapreserving-2}) the length-preserving evolution is faster in the outward normal direction and the limiting curve is approaching a larger circle with the perimeter equal to the length of the initial curve. 
Again, due to application of the curvature adjusted tangential velocity we were able to accurately handle initial large variations in the curvature.

\begin{figure}
\begin{center}
\includegraphics[width=0.52\textwidth]{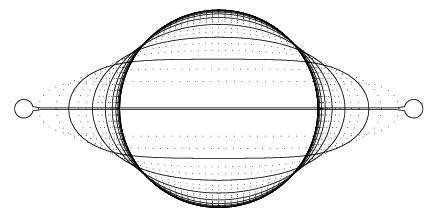}
\end{center}

\caption{Evolution of a thin dumb-bell initial curve by the isoperimetric ratio gradient flow.}
\label{fig:isogradient}
\end{figure}

\begin{figure}
\begin{center}
\includegraphics[width=0.4\textwidth]{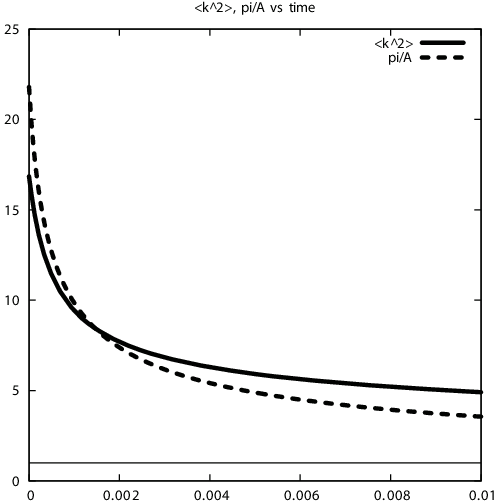}
\end{center}

\caption{Comparison of $\langle k^2 \rangle = \frac1L \int_\Gamma k^2\ds$ (solid line) and $\frac{\pi}{A}$ (dotted line). Evolution of the thin dumb-bell initial curve by the isoperimetric ratio gradient flow. Initial violation of the Gage inequality due to nonconvexity of $\Gamma^t$.}
\label{fig:isogradient-gageinequality}
\end{figure}

In our last set of examples we present the isoperimetric ratio gradient flow. In this flow, the normal velocity is given by 
 \[
\beta= k - \frac{L}{2 A}.
\]
In Fig.~\ref{fig:isogradient} we show evolution of the same thin dumb-bell initial curve as in Figs.~\ref{fig:areapreserving-2},\ref{fig:lengthpreserving-2}. The limiting curve is again a circle. It is a larger (smaller) circle when compared to the area- (length-) preserving flow. The initial thin dumb-bell initial curve is not convex and it violates the Gage isoperimetric inequality (\ref{gageineq}). If we restate inequality (\ref{gageineq}) in the form 
\[
\frac{\pi}{A} \le \frac1L \int_\Gamma k^2 \,\ds \equiv \langle k^2 \rangle
\]
we can observe initial temporal violation of this inequality as it is presented in Fig.~\ref{fig:isogradient-gageinequality}. 

\section*{Conclusions}
In this paper we proposed and analyzed the so-called curvature adjusted tangential velocity for a flow of plane curves. Evolution in the inner normal direction is driven by the normal velocity which may depend on the curvature, position, tangential angle and some nonlocal quantities like the total length, enclosed area and total elastic energy of a curve. We showed local existence, uniqueness and continuation of classical solutions to the system of governing geometric equations. A stable numerical approximation scheme based on the flowing finite volume method with curvature adjusted tangential velocity was also proposed. Its capability has been tested on several computational examples involving nonlocal geometric flows.

\section*{Acknowledgments}
The authors were supported by APVV-0184-10 grant (DS) and 
Grant-in-Aid for Scientific Research (C) 23540150 (SY).

\end{document}